\documentclass[12pt,reqno]{amsart}

\usepackage{amssymb}
\usepackage{bbm}
\usepackage{mathrsfs}
\usepackage{xypic}
\usepackage{float}

\newtheorem{fac}{Fact}[section]
\newtheorem{lem}[fac]{Lemma}
\newtheorem{prop}[fac]{Proposition}
\newtheorem{theo}[fac]{Theorem}
\newtheorem{coro}[fac]{Corollary}
\newtheorem{facs}[fac]{Facts}

\theoremstyle{definition}
\newtheorem{ttt}[fac]{}
\newtheorem{defi}[fac]{Definition}

\theoremstyle{remark}
\newtheorem{rem}[fac]{Remark}
\newtheorem{rems}[fac]{Remarks}
\newtheorem{ex}[fac]{Example}

\newcommand{\pmodulo}[1]{\nobreak\mkern8mu (\textup{mod}\,\,#1)}

\newcommand{\br}{ }
\newcommand{\brr}{, }

\newcommand{\Hom}{\mathop{\text{\rm Hom}}\nolimits}
\newcommand{\Aut}{\mathop{\text{\rm Aut}}\nolimits}
\newcommand{\Div}{\mathop{\text{\rm Div}}\nolimits}
\newcommand{\Pic}{\mathop{\text{\rm Pic}}\nolimits}
\newcommand{\Spec}{\mathop{\text{\rm Spec}}\nolimits}
\newcommand{\Gal}{\mathop{\text{\rm Gal}}\nolimits}
\newcommand{\Br}{\mathop{\text{\rm Br}}\nolimits}
\newcommand{\GL}{\mathop{\text{\rm GL}}\nolimits}

\renewcommand{\div}{\mathop{\text{\rm div}}\nolimits}
\newcommand{\val}{\mathop{\text{\rm val}}\nolimits}
\newcommand{\ram}{\mathop{\text{\rm ram}}\nolimits}
\newcommand{\supp}{\mathop{\text{\rm supp}}\nolimits}

\newcommand{\notd}{\,\mathord{\nmid}\,}

\newcommand{\ev}{\mathop{\text{\rm ev}}\nolimits}
\newcommand{\im}{\mathop{\text{\rm im}}\nolimits}

\newcommand{\tors}{\text{\rm tors}}
\newcommand{\et}{\text{\rm \'et}}

\newcommand{\bbA}{{\mathbbm A}}

\newcommand{\bbD}{{\mathbbm D}}
\newcommand{\bbF}{{\mathbbm F}}
\newcommand{\bbG}{{\mathbbm G}}
\newcommand{\bbN}{{\mathbbm N}}
\newcommand{\bbQ}{{\mathbbm Q}}
\newcommand{\bbR}{{\mathbbm R}}
\newcommand{\bbZ}{{\mathbbm Z}}

\newcommand{\bbx}{{\mathbbm x}}

\newcommand{\calC}{{\mathscr{C}}}

\newcommand{\calE}{{\mathscr{E}}}
\newcommand{\calH}{{\mathscr{H}}}

\newcommand{\calK}{{\mathscr{K}}}
\newcommand{\calL}{{\mathscr{L}}}
\newcommand{\calO}{{\mathscr{O}}}
\newcommand{\calS}{{\mathscr{S}}}
\newcommand{\calU}{{\mathscr{U}}}
\newcommand{\calV}{{\mathscr{V}}}
\newcommand{\calX}{{\mathscr{X}}}

\newcommand{\frakp}{\mathfrak{p}}

\newcommand{\Pb}{{\text{\bf P}}}
\newcommand{\Ab}{{\text{\bf A}}}

\newcounter{abc}
\newenvironment{abc}{\begin{list}{\rm \alph{abc}) }%
{\usecounter{abc} \leftmargin=0.0pt \labelsep=0.0pt %
\listparindent=0.0pt \labelwidth=0.0pt \parsep=\smallskipamount %
\itemsep=0.0pt \topsep=0.0pt \partopsep=\smallskipamount}}{\end{list}}

\newcounter{iii}
\newenvironment{iii}{\begin{list}{\rm \roman{iii}) }%
{\usecounter{iii} \leftmargin=0.0pt \labelsep=0.0pt %
\listparindent=0.0pt \labelwidth=0.0pt \parsep=\smallskipamount%
 \itemsep=0.0pt \topsep=0.0pt \partopsep=\smallskipamount}}{\end{list}}

\def\rightend#1#2{{%
 \leavevmode\nobreak\hskip .5em plus 1fil
 \penalty600 \hskip 0pt plus -1filll
 \vadjust{}\nobreak\hskip 0pt plus 1filll%
 #1\parfillskip=#2\relax \par}}

\def\eop{\ifmmode\rule[-22pt]{0pt}{1pt}\ifinner\tag*{$\square$}\else\eqno{\square}\fi\else\rightend{$\square$}{0pt}\fi}

\newcommand{\ratarrow}{$%
$\definemorphism{rat}\dashed\tip\notip%
\spreaddiagramcolumns{-12pt}%
 - \!\!\diagram%
\rrat & 
\enddiagram\!\!$%
$}

\title[On integral points on degree four del Pezzo surfaces]{On integral points on degree four\\ del Pezzo surfaces}

\begin{document}

\author[J\"org Jahnel]{J\"org Jahnel${}^\ddagger$}

\address{D\'epartement Mathematik\\ Univ.\ \!Siegen\\ \!Walter-Flex-Str.\ \!3\\ D-57068 \!Siegen\\ \!\mbox{Germany}}
\email{jahnel@mathematik.uni-siegen.de}
\urladdr{http://www.uni-math.gwdg.de/jahnel}

\renewcommand{\thefootnote}{\fnsymbol{footnote}}
\author[Damaris Schindler]{Damaris Schindler${}^{*,\ddagger}$}

\address{Mathematisch Instituut\\ Universiteit Utrecht\\ Budapestlaan~6\\ NL-3584 CD Ut\-recht\\ The Netherlands}
\email{d.schindler@uu.nl}
\urladdr{http://www.uu.nl/staff/DSchindler}



\date{\today}



\begin{abstract}
We report on our investigations concerning algebraic and transcendental Brauer-Manin obstructions to integral points on complements of a hyperplane section in degree four del Pezzo surfaces. We~discuss moreover two concepts of an obstruction at an archimedean~place. Concrete examples are given of pairs of non-homogeneous quadratic polynomials in four variables representing
$(0,0)$
over~$\bbQ$
and over
$\bbZ_p$
for all
primes~$p$,
but not
over~$\bbZ$.
By~blow-up, these yield cubic polynomials in three variables all integral solutions of which satisfy a gcd~condition.
\end{abstract}


\footnotetext[1]{The second author was supported by the {\em NSF\/} under agreement No.\ {\tt DMS-1128155}.}

\footnotetext[3]{All computations are with {\tt magma}~\cite{BCP}.}

\maketitle
\thispagestyle{empty}

\section{Introduction}

One says that a variety over a number field
$k$
satisfies the Hasse principle if the existence of a rational point over every
completion~$k_\nu$
of~$k$
is enough to imply the existence of a
\mbox{$k$-rational}~point.
This principle takes its name from the classical Hasse-Minkowski theorem, which states that the Hasse principle holds for the class of quadric hypersurfaces. Already~for cubic curves and cubic surfaces however, the Hasse principle can fail. There has been much work on constructing and controlling such failures, particularly in recent~time.

Many~but not all (cf.~\cite[Theorem~2]{Sk99}) examples of failures of the Hasse principle that are known today can be explained by one and the same argument, which was discovered by Yu.\,I.~Manin~\cite{Ma} around~1970. In~fact, a global Brauer class
$\alpha \in \Br(X)$
is usually responsible for the non-existence of
$k$-rational
points on the underlying proper
variety~$X$
over~$k$.
This~phenomenon is called the Brauer-Manin obstruction to rational~points. Its~mechanism works as~follows.

Let
$\frakp$
be an arbitrary prime ideal
of~$k$.
The~Grothendieck-Brauer group is a contravariant functor from the category of~schemes to the category of abelian groups. In~particular, given a scheme
$X$
and a
\mbox{$k_{\nu_\frakp}$-rational}
point
$x \colon \Spec k_{\nu_\frakp} \to X$,
there is a restriction~homomorphism
$x^*\colon \Br(X) \to \Br(k_{\nu_\frakp}) \cong \bbQ/\bbZ$.
For~a Brauer class
$\alpha \in \Br(X)$,
one calls
$$\ev_{\alpha,\nu_\frakp}\colon X(k_{\nu_\frakp}) \longrightarrow \bbQ/\bbZ \, , \quad x \mapsto x^*(\alpha) \, ,$$
the local evaluation map, associated
with~$\alpha$.
Analogously, when
$\nu$
is a real place, there is the local evaluation map
$\ev_{\alpha,\nu} \colon X(k_\nu) \to \frac12\bbZ/\bbZ$
while, for a complex place, one has the zero~map.

If~$X$
is a proper
\mbox{$k$-scheme}
then the local evaluation maps are continuous with respect to the
\mbox{$\nu$-adic}
topologies
on~$X(k_\nu)$.
Moreover,~$\ev_{\alpha,\nu}$
is constant for all but finitely many places~\cite[Chapter~IV, Proposition~2.3.a.ii)]{Ja}.
Thus,~only adelic points
$x = (x_\nu)_{\nu\in \Omega} \in X(\bbA_k)$
satisfying
$\sum_{\nu\in \Omega} \ev_{\alpha,\nu} (x_\nu) = 0 \in \bbQ/\bbZ$
may possibly be ap\-prox\-imated by
$k$-rational
points.\medskip

\noindent
{\em Brauer-Manin obstruction to integral points.}\smallskip

\noindent
The case of a non-proper
variety~$U$
has been studied only much more recently. Here,~for
$\alpha \in \Br(U)$,
there is usually no finite set
$S$
of places such that
$\ev_{\alpha,\nu}$
would be constant for all
$\nu \not\in S$.
Thus,~the local evaluation maps seem to be useless.

The picture changes however when integral points are considered instead. So~let
$U$
be quasi-projective
over~$k$
and
$\calU$
be a model
of~$U$
that is defined
over~$\Spec \calO_k$,
for
$\calO_k$
the ring of integers
in~$k$.
Then,~similarly to the above, there exists a finite set
$S = S_{\calU,\alpha}$~\cite[\S5.2]{Sk01}
of places, including all archimedean ones, such that the restriction
$$\ev_{\alpha,\nu} |_{\calU(\calO_{k,\nu})} \colon \calU(\calO_{k,\nu}) \longrightarrow \bbQ/\bbZ$$
is the zero map for every place
$\nu \not\in S$.
Consequently, for the set of all
\mbox{$S$-integral}
points, one has the inclusions (cf.\ \mbox{\cite[\S1]{CX}})
\begin{align*}
\calU(\calO_{k,S}) \subset \Big( &\prod_{\nu \in S} U(k_\nu) \times \prod_{\nu \not\in S} \calU(\calO_{k,\nu}) \Big)^{\Br(U)} \subset \prod_{\nu \in S} U(k_\nu) \times \prod_{\nu \not\in S} \calU(\calO_{k,\nu}) \\
\intertext{for}
\Big( &\prod_{\nu \in S} U(k_\nu) \times \prod_{\nu \not\in S} \calU(\calO_{k,\nu}) \Big)^{\Br(U)} := \\
\Big\{ (x_\nu)_\nu \in &\prod_{\nu \in S} U(k_\nu) \times \prod_{\nu \not\in S} \calU(\calO_{k,\nu}) \;\Big|\; \forall \alpha \in \Br(U)\colon \sum_\nu \ev_{\alpha,\nu}(x_\nu) = 0 \Big\} \,.
\end{align*}
In very much the same way as Yu.\,I.~Manin did in his book~\cite{Ma} for counterexamples to the Hasse principle for rational points, J.-L.\ Colliot-Th\'el\`ene and F.\ Xu~\cite{CX} explained many classical counterexamples to the integral Hasse principle or strong approximation (cf.~\cite[Definition~2.1]{DW} or~\cite[\S7.1]{PR}) off certain primes by the Brauer-Manin obstruction to integral points. In \cite{Xu}, F.\ Xu even succeded in showing that the Brauer-Manin obstruction is the only obstruction to strong approximation off infinity for certain quadric fibrations, given that the set of real points is non-compact.

Colliot-Th\'el\`ene and Xu already considered not only algebraic, but also transcendental Brauer classes in their work (see Remark 2.11 and Remark 3.8 in \cite{CX}), and more examples have been constructed recently by A.~Kresch and Yu.~Tschinkel~\cite{KT}.\bigskip\pagebreak[3]

\noindent
{\em Obstruction at infinity.}\smallskip

\noindent
When asking for solutions in integers, there is another point to be considered, besides local solubility and the Brauer-Manin obstruction.
Consider,~for example, the equation
$X^2+23Y^2=2$,
which is insoluble in
$\bbZ^2$,
although there exist rational solutions, as well as
\mbox{$p$-adic}
solutions for every
prime~$p$.
The~most natural argument to show insolubility
over~$\bbZ$
is certainly to notice that the equation defines a compact submanifold
of~$\bbR^2$,
which, as integral points are sought for, leaves only finitely many~cases. The~same effect occurs for every system of equations in an arbitrary number of variables that includes
$Q(X_1,\ldots,X_m) = C$,
for~$Q$
a definite quadratic~form. Cf.~\cite[Exemple~5.9]{CW}, which is exactly of this~type. Such~examples are obstructed at the real place. And~they are in fact obstructed in the crudest possible~sense.

On~the other hand, consider the equation
$X^2-Y^2=2$.
Here,~the submanifold defined by the equation is non-compact, but still there is an obstruction at the real place, in the sense~below. The~equation implies that
$\smash{|X-Y| \leq \sqrt{2}}$
or
$\smash{|X+Y| \leq \sqrt{2}}$.
Thus,~considering only integral points, one is reduced to finitely many algebraic subvarieties of lower~dimension. The~same effect occurs, of course, whenever an equation of the form
$F_1(X_1,\ldots,X_m)F_2(X_1,\ldots,X_m) = C$
is considered, for
$F_1$
and~$F_2$
coprime, non-constant forms in an arbitrary number of variables.

When working over a number field, the latter type of obstruction may even occur at a complex~place. For~example, every solution
in~$\bbZ[i]^2$
of the equation
$X^2+Y^2 = C$
fulfils either
$\smash{|X-iY| \leq \sqrt{|C|}}$
or~$\smash{|X+iY| \leq \sqrt{|C|}}$.

Concerning the concept of unobstructedness at infinity, it is our impression that a definite form for it has yet to be~found. We~think, however, that the examples just considered, despite being so simple, are typical in a certain~sense. Thus,~in section~\ref{sec_obs}, we suggest two notions, called weak and strong unobstructedness, that reflect the effects just indicated, and compare them to results existing in the~literature. We~also show that the assumption of strong unobstructedness is in some sense necessary, when one asks for strong approximation or Zariski-density of integral points.\medskip

\noindent
{\em Our results.}\smallskip

\noindent
The goal of this article is to present some theory and more examples, in which strong approximation off
$\{\infty\}$,
strong approximation off a larger set of primes, or even the integral Hasse principle are violated. We~consider open subschemes of del Pezzo surfaces of degree four that are complements of an effective divisor from the anticanonical class. These~are, in particular, affine
$\log$
$K3$
surfaces of Picard
rank~$5$.
Contrary~to the examples previously studied, in this case both algebraic and transcendental Brauer classes may well~occur.

The article is organised as~follows. Below~this introduction, we first deal with obstructions at infinity (Section~\ref{sec_obs}), followed by generalities about degree four del Pezzo surfaces and complements of hyperplane sections through them (Section~\ref{sec_dP4}), as well as the Brauer-Manin obstruction to integral points (Section~\ref{sec_Brauer}). We~compute the list of all groups that may occur as the algebraic part of the Brauer groups of the surfaces considered in Theorem~\ref{algBr}. This~turns out to comprise
$0$,
$\bbZ/2\bbZ$,
$(\bbZ/2\bbZ)^2$,
$(\bbZ/2\bbZ)^3$,
$(\bbZ/2\bbZ)^4$,
$\bbZ/4\bbZ$,
$\bbZ/2\bbZ \!\times\! \bbZ/4\bbZ$,
and\/
$(\bbZ/2\bbZ)^2 \!\times\! \bbZ/4\bbZ$,
which is much richer than the corresponding list in the case of a proper degree four del Pezzo surface. For~our results on the transcendental part of the Brauer group, see Theorem~\ref{trBr} and Corollary~\ref{trBrcor}.

Section~\ref{sec_explBr} then deals with the question how to evaluate a Brauer class practically. We~present model cases, in which this may be done in  satisfactory manner. Finally, Sections~\ref{sec_ex1}, \ref{sec_ex2}, and \ref{sec_ex3} are devoted to explicit~examples.

In a classical language, these may be formulated in the way~below. A~pair of (non-homogeneous) quadratic polynomials
$q_1, q_2 \in \bbZ[X_0, \ldots, X_3]$
in four variables does not represent
$(0,0)$
in integers, although the system of equations
\begin{align*}
q_1(X_0, \ldots, X_3) &= 0 \\
q_2(X_0, \ldots, X_3) &= 0
\end{align*}
is soluble over
$\bbQ$,
as well as over
$\bbZ_p$
for every prime
number~$p$.

For~instance, in Section~\ref{sec_ex1}, we consider the following system of~equations,
\begin{align*}
           X_0X_1+X_2^2 &= X_3 \, , \\
X_3(2X_1+X_2+X_3)+X_0^2 &= X_1 \, .
\end{align*}
We~show that every
\mbox{$\smash{\bbZ[\frac1{17}]}$-valued}
solution of this system fulfils the relation
$\smash{(x_1, x_3)_2 = 1}$,
for
$(.\,,.)_2$
the
\mbox{$2$-adic}
Hilbert symbol, although the analogous relation is not true for all
\mbox{$\bbZ_2$-valued}
solutions. In~particular, strong approximation 
off~$\{17,\infty\}$
is~violated.

From~this system, it is easy to deduce another one that has no integral solution at all, thereby violating the integral Hasse~principle. Moreover,~blowing up the point
$(0\!:\!1\!:\!0\!:\!0\!:\!0)$,
we are led to the affine cubic surface
\begin{eqnarray*}
128Y_0^3 + 144Y_0^2 + 32Y_0Y_1Y_2 + 8Y_0Y_1 + 128Y_0Y_2^2 + 80Y_0Y_2 + 66Y_0 & \\
{} - 16Y_1^2Y_2- 3Y_1^2 - 4Y_1Y_2 + 80Y_2^2 &+ 40Y_2 + 12 = 0
\end{eqnarray*}
with the amusing property that each of its integral points must satisfy the condition that
$\gcd(8y_0+3, 16y_2+3) > 1$.

Neither~of these examples can be explained by an obstruction at the real place, even in a strong~sense. Observe,~for instance, that the homogeneous parts
$\smash{\underline{q}_1, \underline{q}_2}$
of
degree~$2$
always generate a pencil consisting entirely of indefinite quadratic~forms. This~implies that the real manifold defined by
$(q_1, q_2)$
is non-compact \cite[Satz~1]{Fi} (see also \cite[Tables 1--3]{TWMW}) and, hence, has limit points at the~boundary. In~addition to that, both our surfaces in Sections~\ref{sec_ex1} and~\ref{sec_ex2} are strongly unobstructed
at~$\infty$.
The~arithmetic behaviour we describe is entirely explained by algebraic and transcendental Brauer classes. In~this respect, our examples are similar to the ones given in~\cite{CX}, but different from those recently presented by Y.~Harpaz~\cite[\S4.2]{Harp} and U.~Derenthal and D.~Wei~\cite[Example 6.2]{DW}.

The surface we examine in Section~\ref{sec_ex3} is somewhat~different. Its~arithmetic is explained by Brauer classes together with effects caused by the underlying real
manifold~$U(\bbR)$
being disconnected into non-compact components and a compact~one. So,~similarly to \cite[Example~6.2]{DW}, an interaction of effects caused by Brauer classes occurs with an obstruction at~infinity.
This~example shows that any reasonable definition of unobstructedness at infinity must include requirements on all connected components. In~particular, it is insufficient just to require the existence of a single well-behaved~one.\medskip

\noindent
{\bf Acknowledgements:} We thank Prof.\ J.-L.\ Colliot-Th\'el\`ene for comments on an earlier version of this paper and Dr.\ M.\ Bright for pointing out an oversight in our original formulation of Corollary \ref{trBrcor}. The second author is supported by a NWO grant 016.Veni.173.016.

\section{Obstruction at infinity}
\label{sec_obs}

\begin{ex}
\label{obst}
Consider the subscheme
$\calX \subset \Pb^4_\bbZ$
defined by the system of~equations
\begin{align*}
2X_0^2+X_1^2+X_2^2 &= 26X_4^2 \, , \\
3X_1^2+X_2^2+X_3^2 &= 13X_4^2
\end{align*}
and put
$\calU := \calX \!\setminus\! \calH$,
for
$\calH := V(X_4) \subset \Pb^4_\bbZ$.
Then~the generic fibre
$X$
of~$\calX$
is a degree four del Pezzo surface.
Moreover,~there exist
\mbox{$\bbQ$-rational}
points on the generic fibre
$U$
of~$\calU$,
as well as
\mbox{$\bbZ_p$-valued}
points
on~$\calU$
for every prime
number~$p$.
Examples are
$\smash{(\frac{18}7\!:\!\frac{1}7\!:\!\frac{25}7\!:\!\frac{3}7\!:\!1)}$
and~$\smash{(\frac{54}{19}\!:\!\frac{23}{19}\!:\!\frac{55}{19}\!:\!\frac{9}{19}\!:\!1)}$.
Yet,~there is no integral~point.
\end{ex}

This~is, however, not at all an interesting~example. The reason is simply that there is an obstruction at infinity. Indeed,~the two equations imply
$|x_0|, \ldots, |x_3| \leq 3$
for every
\mbox{$\bbQ$-rational}
point
$(x_0\!:\!\ldots\!:\!x_3\!:\!1)$
on~$\calU$.
For~integral points, this leaves us with only finitely many cases and it turns out that actually none of them corresponds to an integral~point. This motivates the first part of the following definition. The~second part is motivated by examples of the type
$X^2-Y^2=2$,
as discussed in the~introduction.

\begin{defi}
\label{def_obst}
Let~$k$
be a number field and
$\nu$
an archimedean place
of~$k$.
Furthermore,~suppose that
$X \subset \Pb^n_k$
is an irreducible, closed subscheme,
$l \in \Gamma(\Pb^n_k, \calO(1))$
a linear form, and
$H := V(l) \subset \Pb^n_k$
the corresponding
\mbox{$k$-rational}
hyperplane. Put
$U := X \!\setminus\! H$.

\begin{iii}
\item
We say that
$U$
is {\em weakly unobstructed at
$\nu$\/}
if every connected component of
$U(k_\nu)$
has a limit point
on~$H$.
\item
We~call
$U$
{\em strongly unobstructed at
$\nu$\/}
if every connected component
$U'$
of~$U(k_\nu)$
fulfils the following~condition.

Given an integer
$d > 0$,
finitely many
\mbox{$k$-rational}
forms
$s_1, \ldots, s_r \in \Gamma(X, \calO(d))$,
neither being a scalar multiple
of~$l^d$,
and a constant
$c \in \bbR$,
there always exists a point
$x \in U'$
such that
$\| s_i(x)/l^d(x) \|_\nu \geq c$
for
all~$i \in \{1, \ldots, r\}$.
\item 
If~$U$
is not weakly unobstructed (strongly unobstructed) at
$\nu$
then we say that
$U$
is {\em strongly obstructed (weakly obstructed)\/}
at~$\nu$.
\end{iii}
\end{defi}

\begin{rems}
\begin{iii} 
\item
As~$\Pb^n(k_\nu)$
is a compact topological space, the property of being weakly unobstructed simply means that
$U(k_\nu)$
does not have a bounded connected~component. We~note that if
$\nu$
is a complex place and
$\dim U \geq 1$
then
$U$
is weakly unobstructed
at~$\nu$.
\item
On~the other hand, strong unobstructedness states that no connected component of
$U(k_\nu)$
is contained in a finite union of tubular neighbourhoods of hypersurfaces of the form
$\| s_i(x)/l^d(x) \|_\nu \leq c$.
\item
Example~\ref{obst} is strongly obstructed. In~fact, the topological boundary of
$\calU(\bbR)$
is given by
$2X_0^2+X_1^2+X_2^2 = 3X_1^2+X_2^2+X_3^2 = 0$,
which defines the empty set
in~$\Pb^3(\bbR)$.
\item
As~already remarked in the introduction, it seems to us that a definite form for the concept of unobstructedness at infinity has yet to be~found. For~instance, in~\cite[Theorems~1.1 and~1.2]{DW}, U.~Derenthal and D.~Wei included some unboundedness conditions at infinity and showed that they are, in general,~necessary. In the case of a fibration over the affine
line~$\Ab^1$,
our notion of being strongly unobstructed implies unboundedness in the sense of Derenthal and~Wei. Our~definition is partly motivated by their examples, which however come with a natural projection to some affine~space. Our definition aims to suggest a general notion of unobstructedness at infinity, which does not require such a fibration.

It~is inspired, too, by the results of J.-L.\ Colliot-Th\'el\`ene and F.~Xu~\cite[Theorems~3.7, 4.5, and~8.3]{CX}, as well as by the recent work of Y.~Harpaz~\cite{Harp}. In~fact, in the case that
$U(k_\nu)$
is connected, Harpaz' concept of being split~\cite[Def.~1.0.8]{Harp} is~slightly stronger than weak unobstructedness. It~seems, however, that split varieties do not behave well, in general, when
$U(k_\nu)$
is~disconnected, cf.~Section~\ref{sec_ex3}. Y.~Harpaz' example of a surface that is split, has no Brauer-Manin obstruction, yet no integral point, is weakly obstructed at infinity, according to our definition, as we show~below. It is exactly this property of being obstructed that is used in his work to show that the surface has no integral~points.
\end{iii}
\end{rems}

\begin{ex}[cf.~{\cite[\S4.2]{Harp}}]\label{Harpazex}
Let
$X \subset \Pb^3_\bbQ$
be defined by the equation
$$
((11X_0+5X_3)X_1+3X_3^2)X_2 = (3X_0+X_3)X_3^2
$$
and put
$U := X \!\setminus\! H$,
for
$H := V(X_3) \subset \Pb^4_\bbQ$.
Then~$U$
is weakly unobstructed at
$\infty$,
but not strongly~unobstructed.\smallskip

\noindent
{\bf Proof.}
It is not hard to see that
$U(\bbR)$
is connected. Hence,~weak unobstructedness follows from
$U$
being split in the terminology of~\cite[Def.~1.0.8]{Harp}.

In~order to disprove strong unobstructedness, observe that every real point
$x = (x_0\!:\!x_1\!:\!x_2\!:\!1) \in U(\bbR)$
must fulfil the equation
$((11x_0+5)x_1+3)x_2 = 3x_0+1$,
i.e.\
$\smash{x_0(11-\frac3{x_1x_2}) = \frac1{x_1x_2} -\frac3{x_1} -5}$.
This immediately shows that
$|x_1|, |x_2| \geq 1$
implies
$|x_0| \leq \frac98$.
\eop
\end{ex}

\begin{ex}
Let
$X \subset \Pb^6_\bbQ$
be defined by the equation
$$
(X_0^2 + X_1^2 - X_2^2)(X_3^2 + X_4^2 - X_5^2) = (X_0^2 + X_1^2 - X_2^2 + X_3^2 + X_4^2 - X_5^2)X_6^2
$$
and put
$U := X \!\setminus\! H$,
for
$H := V(X_6) \subset \Pb^4_\bbQ$.
Then~$U$
is weakly unobstructed at
$\infty$,
but not strongly~unobstructed.\smallskip

\noindent
{\bf Proof.}
It~is not hard to see that
$U(\bbR)$
has exactly two connected components and that
$(1\!:\!0\!:\!0\!:\!0\!:\!0\!:\!0\!:\!0)$
is a limit point for either of~them.

On~the other hand, every real point
$x = (x_0\!:\!x_1\!:\!x_2\!:\!x_3\!:\!x_4\!:\!x_5\!:\!1) \in U(\bbR)$
must fulfil the equation
$Q(x_0, x_1, x_2)Q(x_3, x_4, x_5) = Q(x_0, x_1, x_2) + Q(x_3, x_4, x_5)$,
and hence
$|Q(x_0, x_1, x_2)| \leq 2$
or
$|Q(x_3, x_4, x_5)| \leq 2$
for~$Q$
the ternary quadratic form
$Q(X_0,X_1,X_2) := X_0^2 + X_1^2 - X_2^2$.
\eop
\end{ex}

The next theorem shows that the notion of being strongly unobstructed is essentially necessary. In fact otherwise, under the assumption that $U(\bbR)$ is connected, the integral points
on~$\calU$
are contained in a union of finitely many lower-dimensional subschemes. In~particular, they are not Zariski dense and strong approximation
off~$\{\infty\}$
cannot be~fulfilled.

\begin{theo}
\label{weaklyobstr}
Let\/~$\calX \subset \Pb^n_\bbZ$
be an irreducible, closed subscheme,
$l \in \Gamma (\Pb^n_\bbZ, \calO(1))$
a~linear form,
$\calH := V(l) \subset \Pb^n_\bbZ$,
and\/~$\calU:= \calX\setminus \calH$.
Denote~the generic fibres
of\/
$\calX$
and\/~$\calU$
by\/
$X$
and\/~$U$,
respectively.\smallskip

\noindent
Assume that\/
$U(\bbR)$
is connected and that\/
$U$
is weakly obstructed at~infinity. Then the set\/
$\calU(\bbZ)$
of integral points is not Zariski-dense
in\/~$\calU$.\medskip

\noindent
{\bf Proof.}
{\em
Since
$U$
is weakly obstructed at infinity, there exist an integer
$d>0$
and forms
$s_1, \ldots, s_r \in \Gamma (X, \calO_X(d))$,
none of them a multiple of
$l^d$,
as well as a positive constant
$c>0$,
such that the following is~true. For~every
$x \in U(\bbQ)$,
there is an
$1 \leq i\leq r$
such that
$$\left| \frac{s_i(x)}{l^d(x)} \right| \leq c \, .$$
Taking global sections commutes with arbitrary flat base extensions. In~particular,
$\Gamma(X, \calO_X(d)) = \Gamma(\calX\!, \calO_\calX(d)) \otimes_\bbZ \bbQ$.
As~a consequence of this, we see that there exists a non-zero integer
$M \in \bbZ$
such that
$Ms_1, \ldots, Ms_r \in \Gamma (\calX\!, \calO_\calX(d))$.

Consider now an integral point
$x \in \calU(\bbZ)$.
Its restriction to the generic fibre has the property that
$$\left| \frac{Ms_i(x)}{l^d(x)} \right| \leq Mc \, ,$$
for some
$1 \leq i \leq r$.
Thus,~for every such integral point, there exists an integer
$B$
such that
$|B| \leq Mc$
and
$$\frac{Ms_i(x)}{l^d(x)} = B \, .$$
In particular, this system defines a Zariski-closed subset
$Z \subset U$.
Moreover, if we had
$Z=U$
then
$s_i$
would be equal to
$B$
on
$U$,
a contradiction to none of the
$s_i$
being a scalar multiple of the
form~$l^d$.
Hence,
$Z \subsetneqq U$.

We conclude that there is a finite union of such proper Zariski-closed subsets of~$U$
that contains all restrictions of integral points
$P \in \calU(\bbZ)$
to the generic fibre. Hence, the integral points on
$\calU$
cannot be Zariski-dense.
}
\eop
\end{theo}

We note that Theorem \ref{weaklyobstr} may be seen as a generalisation of Theorem 4.1.1. of Harpaz's work \cite{Harp}. Finally, we close this section in giving a useful criterion to detect varieties that are strongly unobstructed. 

\begin{theo}
\label{unobstr}
Let\/~$k$
be a number field,
$\nu$~an
archimedean place
of\/~$k$,
$X \subset \Pb^n_k$
be a normal, projective variety,
$l \in \Gamma(\Pb^n_k, \calO(1))$
a linear form,
$H := V(l) \subset \Pb^n_k$
the corresponding\/
\mbox{$k$-rational}
hyperplane, and put\/
$U := X \!\setminus\! H$.\smallskip

\noindent
Suppose~that

\begin{iii}
\item[$\bullet$ ]
the scheme\/
$(H \cap X)_{k_\nu}$
is reduced and irreducible and that
\item[$\bullet$ ]
every connected component
of\/~$U(k_\nu)$
has a limit point\/
$x \in (H \cap X)(k_\nu)$
that is non-singular as a point
on\/~$H \cap X$.
\end{iii}

\noindent
Then\/~$U$
is strongly unobstructed
at\/~$\nu$.
\end{theo}

Note~that Theorem~\ref{unobstr} does not apply to example \ref{Harpazex}, despite
$X$
being a normal cubic~surface, since
$H \cap X$
is a union of three~lines.\medskip

\noindent
{\bf Proof of Theorem~\ref{unobstr}.}
If~$\nu$
is a complex place then
$U(k_\nu)$
is clearly connected, but the same is not necessarily true for a real~place. Thus,
let~$U'$
be a connected component
of~$U(k_\nu)$.
There~is a limit point
$x' \in (H \cap X)(k_\nu)$
that is non-singular on the limit set. Therefore,~there exists a neighbourhood
$V \ni x'$
such that
$[U' \cup (H \cap X)(k_\nu)] \cap V$
is diffeomorphic to some
$\bbD^m$,
$\bbR^m$
or half space
$\bbR^{m-1} \!\times\! \bbR_{\geq 0}$.

Now~let
degree~$d$
forms
$s_1, \ldots, s_r$
not being multiples
of~$l^d$
and a constant
$c$
be given. We~choose a linear form, which we denote
by~$l'$,
in such a way that
$l'(x') \neq 0$.
Shrinking~$V$,
if necessary, we may assume that the same holds on the whole
of~$V$.

By~our assumptions,
$(H \cap X)_{k_\nu}$
defines a prime Weil divisor
on~$X_{k_\nu}$.
Put
$d_i := \nu_{(H \cap X)_{k_\nu}} (s_i)$
for~$i = 1, \ldots, r$.
Since~$s_1, \ldots, s_r$
are not multiples
of~$l^d$,
we have
$0 \leq d_i < d$
for
each~$i$.
The~rational functions
$$f_i := \frac{s_i}{l^{d_i}l'^{d-d_i}}$$
have neither zeroes nor poles along the divisor
$(H \cap X)_{k_\nu}$.

In~particular, there exists a point
$x'' \in (H \cap X)(k_\nu) \cap V$
such that
$f_i(x'') \neq 0$
for~$i = 1, \ldots, r$.
Say,~we have
$|s_i(x'')/l^{d_i}(x'')l'^{d-d_i}(x'')| \geq \varepsilon$
for every
$i = 1, \ldots, r$
and
some~$\varepsilon > 0$,
which we may assume to
fulfil~$\varepsilon < 2c$.
Moreover,~$l(x'') = 0$,
since
$x'' \in (H \cap X)$.

Next,~we approximate
$x''$
by a point
$x''' \in U' \cap V$
outside the limit~set. We~may enforce that
$|s_i(x''')/l^{d_i}(x''')l'^{d-d_i}(x''')| \geq \varepsilon/2$
for every
$i = 1, \ldots, r$
and that
$0 < | l(x''')/l'(x''')| \leq \frac\varepsilon{2c}$.
This~yields
$$\left| \frac{s_i(x''')}{l^d(x''')} \right| = \left| \frac{s_i(x''')}{l^{d_i}(x''')l'^{d-d_i(x''')}} \right| \!\cdot\! \left| \frac{l'(x''')}{l(x''')} \right|^{d-d_i} \geq \frac\varepsilon2 \cdot \left( \frac{2c}\varepsilon \right)^{d-d_i} \geq \frac\varepsilon2 \cdot \frac{2c}\varepsilon = c$$
for
$i = 1, \ldots, r$,
as~required.
\eop

\section{Del Pezzo surfaces of degree four and\\ complements of hyperplane sections}
\label{sec_dP4}

We are interested in quasi-projective schemes
$U \subset \Pb^4_k$,
given by systems of equations of the type
\begin{eqnarray*}
Q_1(X_0, \ldots, X_4) &=& 0 \, , \\
Q_2(X_0, \ldots, X_4) &=& 0 \, , \\
X_4 &\neq& 0 \, ,
\end{eqnarray*}
for homogeneous quadrics
$Q_1, Q_2 \in \calO_k[X_0, \ldots, X_4]$.
Then~a rational point
on~$U$
is the same as a solution
$(x_0, \ldots, x_3, 1)$
for
$x_0, \ldots, x_3 \in k$,
while an integral point may be thought of as a solution
$(x_0, \ldots, x_3, 1)$
for
$x_0, \ldots, x_3 \in \calO_k$.
The~scheme
$U$
itself is an open subscheme of a del Pezzo surface of degree~four.

By definition, a del Pezzo surface is a non-singular, proper algebraic
surface~$X$
over a field
$k$
with an ample anti-canonical sheaf
$\calK^{-1}$.
Over~an algebraically closed field, every del Pezzo surface of degree
$d \leq 7$
is isomorphic
to~$\Pb^2$,
blown up in
$(9-d)$
points in general position~\cite[The\-o\-rem~24.4.iii)]{Ma}.

According to the adjunction formula, a non-singular complete intersection of two quadrics in
$\Pb^4$
is del Pezzo. The converse is true, as~well. For every del Pezzo surface of degree four, its anticanonical image is the complete intersection of two quadrics in
$\Pb^4$~\cite[Theorem~8.6.2]{Do}.

Thus,~associated with a degree four del Pezzo
surface~$X$,
there is a pencil
$\smash{(\mu Q_1 + \nu Q_2)_{(\mu:\nu) \in \Pb^1}}$
of quadrics
in~$\Pb^4$,
the base locus of which
is~$X$.
This pencil is uniquely determined up to an isomorphism
of~$\Pb^4$.
It~contains exactly five degenerate fibres, each of which is of rank precisely~$4$~\cite[Proposition~3.26.iv)]{Wi}.

Over an algebraically closed field, a degree four del Pezzo
surface contains exactly 16~lines~\cite[Theorem~26.2.iii)]{Ma} and exactly ten one-dimensional linear systems of conics~\cite[\S2.3]{VAV}. The~latter may be obtained as~follows. Take the five degenerate quadrics in the pencil. As~they are of
rank~$4$,
each of them contains two one-dimensional linear systems of planes. Intersecting~with a second quadric from the pencil, one finds one-dimensional linear systems of conics.\medskip

\begin{facs}[The geometric Picard group]
\label{geomPic}
Let\/~$k$
be an algebraically closed field and\/
$X$
a del Pezzo surface of degree four
over\/~$k$.
Then

\begin{iii}
\item
$\Pic(X) \cong \bbZ^6$.
The~classes of the 16 lines form a generating system
for~$\Pic(X)$.
Moreover,~the canonical class\/
$\calK \in \Pic(X)$
is not divisible by any
integer\/~$\neq \pm1$.
\item
The orthogonal complement\/
$\Pic^0(X) := \langle\calK\rangle^\perp \subset \Pic(X)$
is, up to sign, isometric to the
lattice~$D_5$.
\item
The
group\/~$W$
of all automorphisms
of\/~$\Pic(X)$
respecting\/
$\calK$
and the intersection pairing is isomorphic to the Weyl
group\/~$W(D_5)$.
\item
The operation
of\/~$W$
permutes the classes of the 16\/ lines
in\/~$\Pic(X)$.
This~operation is faithful. The image of the corresponding group homomorphism\/
$W \to S_{16}$
is the largest subgroup
of\/~$S_{16}$
that respects the intersection~matrix.
\item
The operation
of\/~$W$
permutes the classes
in\/~$\Pic(X)$
of the 10 linear systems of~conics. This~operation is faithful,~too. The image of the corresponding group homomorphism\/
$W \to S_{10}$
is isomorphic to\/
$(\bbZ/2\bbZ)^4 \rtimes S_5$.
Here,~$(\bbZ/2\bbZ)^4 \subset (\bbZ/2\bbZ)^5$
means the subgroup
$$\{ (x_0, \ldots, x_4) \in (\bbZ/2\bbZ)^5 \mid x_0 + \ldots + x_4 = 0 \}$$
being acted upon
by\/~$S_5$
in the natural~way.
\end{iii}\smallskip

\noindent
{\bf Proof.}
{\em
i) follows immediately from the description of
$\Pic(X)$
given in~\cite[Proposition~2.2]{VAV}, together with \cite[Theorem~26.2.i)]{Ma}, while ii), iii) and~iv) are parts of~\cite[Theorem~23.9]{Ma}.\smallskip

\noindent
Finally,~v) is known for a long time, as~well. It~follows, for example, from the discussion in \cite[pp.~8--10]{KST}.
}
\eop
\end{facs}

\begin{rem}
The~existence of an isomorphism
$W(D_5) \cong (\bbZ/2\bbZ)^4 \rtimes S_5$
may, of course, be seen directly from the construction of the root system~\cite[\S12]{Hu}. Moreover,
$\Aut(D_5) \cong (\bbZ/2\bbZ)^5 \rtimes S_5$,
in which
$W(D_5)$
is a subgroup of
index~$2$.
\end{rem}

\begin{lem}
\label{lattice}
Let\/~$k$
be an algebraically closed field,
$X$
a del Pezzo surface of degree four
over\/~$k$,
and\/
$D \subset X$
an irreducible divisor such that\/
$\calO(-D) = \calK \in \Pic(X)$
is the canonical~class. Put\/
$U := X \!\setminus\! D$.
Then

\begin{iii}
\item
$\smash{\Pic(U) \cong \Pic(X) / \langle\calK\rangle}$.
\item
The orthogonal projection\/
$\pi\colon \Pic(X) \to \langle\calK\rangle^\perp \otimes_\bbZ \bbR$
induces an~isomorphism
$$\Pic(U) \stackrel{\cong}{\longrightarrow} D_5^* \, ,$$
for\/
$D_5^*$
the lattice dual
to\/~$D_5$.
\end{iii}\smallskip

\noindent
{\bf Proof.}
{\em
i)
This~is \cite[Proposition~II.6.5.c)]{Hart}.\smallskip

\noindent
ii)
By~Fact \ref{geomPic}, the kernel of
$\pi$
is exactly
$\langle\calK\rangle$.
Moreover,~for every
$\calL \in \Pic(X)$,
one has
$\langle \pi(\calL), \calC \rangle = \langle \calL, \calC \rangle \in \bbZ$
for all
$\calC \in \langle\calK\rangle^\perp = D_5$.
Hence~$\pi(\calL) \in D_5^*$
and, consequently,
$\im \pi \subseteq D_5^*$.

In~order to see that in fact equality is true, we first recall that
$D_5$
is of
discriminant~$4$~\cite[\S11, Exc.~2]{Hu}.
Thus,~$\langle\calK\rangle^\perp$
is of
discriminant~$(-4)$.
As~$\langle \calK, \calK \rangle = 4$,
this yields that
$\langle\calK\rangle \oplus \langle\calK\rangle^\perp$
is of
discriminant~$(-16)$.
Moreover,~$\Pic(X)$
is unimodular. Therefore,
$\langle\calK\rangle \oplus \langle\calK\rangle^\perp$
must be of
index~$4$
in~$\Pic(X)$.
Hence~$\# \im\pi/D_5 = 4$.
But~$\# D_5^*/D_5 = 4$,
again by~\cite[\S11, Exc.~2]{Hu}.
}
\eop
\end{lem}

\begin{coro}
Let\/~$X$
and\/~$U$
be as above and\/
$G$
be a group operating on the 16~lines
of\/~$X$
such that the intersection matrix is~respected.

\begin{iii}
\item
Then~the induced operation
of\/~$G$
on\/~$\Pic(X)$
respects the intersection pairing as well as the canonical~class.
\item
The induced operation
of\/~$G$
on\/
$\Pic(U) \cong D_5^*$
takes place via isometries. The~associated homomorphism\/
$G \to \Aut(D_5^*)$
factors
via\/~$W(D_5)$.
\end{iii}\smallskip

\noindent
{\bf Proof.}
{\em
i)
The operation
of~$G$
on the free
\mbox{$\bbZ$-module}
$\Div(X)$
over the 16~lines respects the intersection pairing and, by~Fact~\ref{geomPic}.i),
$\Pic(X)$
is the quotient of
$\Div(X)$
modulo its~radical. Finally,~the canonical class is automatically respected, as the divisor being the sum over all 16~lines defines the invertible sheaf
$\calK^{\otimes (-4)}$.\smallskip

\noindent
ii)
The~first assertion is a direct consequence of~i). The~second follows from Fact \ref{geomPic}.iii).
}
\eop
\end{coro}

Under some minor assumptions,
$U$
is an example of a
$\log$
$K3$
surface. Let~us shortly recall the definition.

\begin{defi}[cf.~{\cite[Def.~2.0.15]{Harp}}]
A
{\em
log
$K3$
surface\/} is a
scheme~$U$
that may be written in the form
$U = X \!\setminus\! D$,
for
$X$
a non-singular, connected scheme that is proper of dimension two over a field
$k$
and
$D \subset X$
an effective divisor such that

\begin{iii}
\item
$\calO(-D) \in \Pic(X)$
is the canonical class and
\item
$\pi_1^\et(U_{\overline{k}}, \;.\;) = 0$.
\end{iii}
\end{defi}

\begin{theo}
\label{logK3}
Let\/~$X \subset \Pb^4$
a non-singular complete intersection of two quadrics over a
field\/~$k$
of
characteristic\/~$0$
and\/
\mbox{$U := X \!\setminus\! H$}
be the complement of a non-sin\-gu\-lar hyperplane~section. Then\/
$U$
is a\/
$\log$
$K3$
surface.\medskip

\noindent
{\bf Proof.}
{\em
We~have
$\pi_1^\et(X_{\overline{k}}, \;.\;) = 0$,
since
$X_{\overline{k}}$
is isomorphic
to~$\Pb^2$,
blown up in five points. Moreover,~\cite[Lemma~3.3.6]{Harp} yields that
$\pi_1^\et(U_{\overline{k}}, \;.\;) = 0$,~too.
Indeed,~each of the 16~lines
on~$X_{\overline{k}}$
intersects the irreducible divisor
$H \cap X$
transversely and in exactly one~point.
Finally,~the adjunction formula shows that
$\calK_X \cong \calO_{\Pb^4}(-5H+2H+2H) |_X \cong \calO_X(-H)$,
which implies the claim.
}
\eop
\end{theo}

\section{The Grothendieck-Brauer group}
\label{sec_Brauer}

\subsubsection*{Generalities.}\leavevmode\medskip

\noindent
By definition~\cite[Remarque~2.7]{GrBrII}, the cohomological Grothendieck-Brauer group
$\Br(U)$
of a scheme
$U$
is the \'etale cohomology group
$H^2_\et(U, \bbG_m)$.
If~$U$
is defined over a field
$k$
then the Hochschild-Serre spectral~sequence~\cite[Exp.\ VIII, Proposition~8.4]{SGA4}
$$H^p(\Gal(\overline{k}/k), H^q_\et (U_{\overline{k}}, \bbG_m)) \Longrightarrow H^{p+q}_\et (U, \bbG_m)$$
carries essential information
about~$\Br(U)$.
It~yields a three-step filtration
$$0 \subseteq \Br_0(U) \subseteq \Br_1(U) \subseteq \Br(U) \, ,$$
the subquotients of which we are going to describe~below.
In~order to avoid unnecessary complications, let us assume that
\begin{equation}
\label{Gm}
\Gamma_\et (U_{\overline{k}}, \bbG_m) = \overline{k} {}^* \, .
\end{equation}

\begin{lem}
\label{algbr}
Let\/~$k$
be a field,
$X$
a non-singular, proper scheme
over\/~$k$,
and\/
$U = X \!\setminus\! D$
be the complement of a
divisor\/~$D$.
Assume~that\/
$D$
splits geometrically into irreducible components\/
$D_1, \ldots, D_r$,
the intersection matrix of which has
rank\/~$r$.\smallskip

\noindent
Then~\/(\ref{Gm}) is~fulfilled.\medskip

\noindent
{\bf Proof.}
{\em
The group
$\Gamma_\et (U_{\overline{k}}, \bbG_m)$
consists of all rational functions
$s$
on~$X_{\overline{k}}$
such that
$\div s = k_1 [D_1] + \ldots + k_r [D_r]$
for some
integers~$k_1, \ldots, k_r$.
Since~$\div s$
is a principal divisor
on~$X_{\overline{k}}$,
its intersection numbers with
$D_1, \ldots, D_r$
must all~vanish. By~our assumption, this is possible only for
$k_1 = \ldots = k_r = 0$.
I.e.,~$s$
must be constant and one sees
that~$\smash{\Gamma_\et (U_{\overline{k}}, \bbG_m) = \overline{k}^*}$,
as~required.
}
\eop
\end{lem}

\begin{rem}
The assumptions on
$D$
are satisfied, in particular, if
$D$
is a geometrically irreducible divisor such that
$D^2 \neq 0$.
\end{rem}

\begin{ttt}
\label{Brsubq}
\begin{iii}
\item
In~general,
$\Br_0(U)$
is the image of a natural homomorphism
\begin{equation}
\label{Leray}
H^2(\Gal(\overline{k}/k), \Gamma_\et (U_{\overline{k}}, \bbG_m)) \longrightarrow \Br(U) \, .
\end{equation}
Under the assumption of~(\ref{Gm}), this shows that
$\Br_0(U)$
is the image of the natural homomorphism
$\Br(k) \to \Br(U)$.
The~homomorphism (\ref{Leray}) is an injection, thereby ensuring that
$\Br_0(U) \cong \Br(k)$,
under some minor assumption. A~sufficient condition is that
$U$
has a
\mbox{$k$-rational}
point or, if
$k$
is a number field, that
$U$
has an adelic point~(cf.~\cite[Proposition~1.3.4.1]{Co}).

The part
$\Br_0(U) \subseteq \Br(U)$
does not contribute to the Brauer-Manin obstruction.
\item
The subquotient
$\Br_1(U)/\Br_0(U)$
is, in general, isomorphic to
$$\ker d_2^{1,1}\colon H^1(\Gal(\overline{k}/k), \Pic(U_{\overline{k}})) \longrightarrow H^3(\Gal(\overline{k}/k), \Gamma_\et (U_{\overline{k}}, \bbG_m)) \, ,$$
for
$d_2^{1,1}$
the differential in the Hochschild-Serre spectral~sequence.

Under~the assumption of~(\ref{Gm}), the right hand side simplifies to
$H^3(\Gal(\overline{k}/k), \overline{k}^*)$.
Moreover,~if
$k$
is a number field then, as 
a by-product of class field theory~\cite[section~11.4]{Ta}, it is known that
$\smash{H^3(\Gal(\overline{k}/k), \overline{k}^*) = 0}$.
Consequently,~one has in~fact
$$\Br_1(U)/\Br_0(U) \cong H^1(\Gal(\overline{k}/k), \Pic(U_{\overline{k}})) \, .$$
The subquotient
$\Br_1(U)/\Br_0(U)$
is called the algebraic part of the Brauer group. The effects of algebraic Brauer classes towards the Brauer-Manin obstruction are called the algebraic Brauer-Manin obstruction.
\item
The subgroup
$\Br_1(U)$
is nothing but the kernel of the natural homomorphism
$\Br(U) \rightarrow \Br(U_{\overline{k}})$.
Thus,~there is a natural injection
$$\Br(U)/\Br_1(U) \hookrightarrow \Br(U_{\overline{k}})^{\Gal(\overline{k}/k)} \, .$$
It~seems to be hard, even in concrete cases, to decide which Galois invariant Brauer classes on
$U_{\overline{k}}$
actually descend
to~$U$.
For~interesting particular results, the reader is advised to study the paper~\cite{CS} of J.-L.\ Colliot-Th\'el\`ene and A.\,N.~Skorobogatov.

The quotient
$\Br(U)/\Br_1(U)$
is called the transcendental part of the Brauer group. The effects of transcendental Brauer classes towards the Brauer-Manin obstruction are called the transcendental Brauer-Manin obstruction.
\end{iii}
\end{ttt}

\subsubsection*{The algebraic part of the Brauer group.}

\begin{ttt}
Let~$k$
be a number field,
$\smash{X \subset \Pb^4_k}$
a degree four del Pezzo surface,
$H \subset \Pb^4_k$
a
\mbox{$k$-rational}
hyperplane such that
$H \cap X$
is geometrically~irreducible,
and~\mbox{$U := X \!\setminus\! H$}.
Then,~$(H \!\cap\! X)^2 = 4$.
Hence,~by Lemma~\ref{algbr} and~\ref{Brsubq}.ii),
$$\Br_1(U)/\Br_0(U) \cong H^1(\Gal(\overline{k}/k), \Pic(U_{\overline{k}})) \, .$$
Moreover,~by Lemma~\ref{lattice}.ii),
$\smash{\Pic(U_{\overline{k}}) \cong D_5^*}$
and, under this isomorphism, the operation
of~$\smash{\Gal(\overline{k}/k)}$
goes over into a homomorphism
$\smash{\Gal(\overline{k}/k) \longrightarrow \Aut(D_5^*)}$
factoring
via~$W(D_5)$.
\end{ttt}

\begin{ttt}
$W(D_5)$~has
197 conjugacy classes of subgroups. It~is known that each of them occurs already for a degree four del Pezzo surface defined
over~$\bbQ$.
Indeed,~this was shown by \smash{B.\,\`E.}~Kunyavskij, A.\,N.\ Skorobogatov, and M.\,A.\ Tsfasman~\cite{KST}. The~reader might want to compare~\cite{EJ}, which reports on a more recent investigation.

We~calculated
$H^1(H, D_5^*)$
for each conjugacy class of subgroups
in~$W(D_5)$,
using {\tt magma}'s functionality for the computation of the cohomology of finite~groups. We~partially reproduce the resulting list below in Figure~1.

\begin{rems}[Some details on the implementation]
\leavevmode\vskip1mm
\begin{iii}
\item
It~should not be difficult to create the
\mbox{$\bbZ W\!(D_5)$-module}
$D_5^*$
within a computer algebra~system. To~be concrete, in {\tt magma}, the command\\[1mm]
{\tt w\_d5 := TransitiveGroup(16,1328);}\\[1mm]
creates the image of
$W(D_5)$
in~$S_{16}$.
Then~the command\\[1mm]
{\tt subs2 := Subsets(Set(GSet(w\_d5)), 2);}\\[1mm]
shows the decomposition of
$\{1,\ldots,16\} \times \{1,\ldots,16\}$
into exactly three orbits under the
\mbox{$W(D_5)$-operation,}
from which the intersection matrix is easily~constructed, cf.~\cite[\S3.1]{JL}. On~the other hand, the command\\[1mm]
{\tt Div := PermutationModule(wd5, Integers());}\\[1mm]
creates a
$\bbZ W\!(D_5)$-module
$\Div$
that, as a
\mbox{$\bbZ$-module},
is free over
$\{1,\ldots,16\}$
and carries the operation
of~$W(D_5)$
that is induced by the one on that~set. Thus,
$D_5^* \cong \Div/\Div_0$,
for
$\Div_0 \subset \Div$
the submodule of rank
$11$
given by the condition that
$\langle x,e_1 \rangle = \ldots = \langle x,e_{16} \rangle$.
Here,
$\langle .,. \rangle$
denotes the pairing encoded by the intersection~matrix.
\item
As our version of {\tt magma} is unwilling to compute the quotient of two
\mbox{$\bbZ W\!(D_5)$-modules}
and allows the user to do this only for
\mbox{$k W\!(D_5)$-modules},
for~$k$
a field, we worked modulo a large prime and recovered the structure of the
\mbox{$\bbZ W\!(D_5)$-module}
$D_5^*$
from this. It turned out that the image of
$W(D_5)$
in
$\GL_5(\bbZ)$
is generated by the two matrices
$$\left(
\begin{array}{rrrrr}
  1 & 0 &-1 & 0 & 1 \\
  0 & 0 & 0 & 1 & 0 \\
 -1 &-1 & 0 &-1 & 0 \\
  0 & 0 & 1 & 0 & 0 \\
 -2 & 0 & 1 &-1 &-1
\end{array}
\right)
\quad{\rm and}\quad
\left(
\begin{array}{rrrrr}
  2 & 1 & \phantom{-}0 & 1 & 1 \\
 -1 &-1 & 0 &-1 & 0 \\
  0 & 1 & 1 & 0 &-1 \\
 -2 & 0 & 1 &-1 &-1 \\
 -1 &-1 & 0 & 0 &-1
\end{array}
\right) .
$$
\end{iii}
\end{rems}\vspace{-0.4cm}

\begin{figure}[H]\vspace{-0.6cm}
\caption{List of algebraic Brauer groups}\vspace{-0.2cm}
\end{figure}
{\tiny
\begin{verbatim}
1, 1, [ 1, 1, 1, 1, 1, 1, 1, 1, 1, 1, 1, 1, 1, 1, 1, 1 ], [],
 []
2, 2, [ 2, 2, 2, 2, 2, 2, 2, 2 ], [ 2 ],
 [ 2, 2, 2 ]
3, 2, [ 2, 2, 2, 2, 2, 2, 2, 2 ], [ 2 ],
 [ 2 ]
4, 2, [ 1, 1, 1, 1, 1, 1, 1, 1, 2, 2, 2, 2 ], [ 2 ],
 []
5, 2, [ 1, 1, 1, 1, 2, 2, 2, 2, 2, 2 ], [ 2 ],
 []
6, 2, [ 2, 2, 2, 2, 2, 2, 2, 2 ], [ 2 ],
 [ 2 ]
7, 3, [ 1, 1, 1, 1, 3, 3, 3, 3 ], [ 3 ],
 []
8, 5, [ 1, 5, 5, 5 ], [ 5 ],
 []
9, 4, [ 4, 4, 4, 4 ], [ 2, 2 ],
 [ 2, 2, 2, 2 ]
10, 4, [ 2, 2, 2, 2, 2, 2, 2, 2 ], [ 2, 2 ],
 []
11, 4, [ 4, 4, 4, 4 ], [ 4 ],
 [ 2, 2, 4 ]
\end{verbatim}

$\vdots$
\begin{verbatim}
181, 96, [ 16 ], [ 2, 2, 2 ],
 [ 2, 2 ]
\end{verbatim}

$\vdots$
\begin{verbatim}
189, 128, [ 16 ], [ 2, 2, 2 ],
 [ 2 ]
190, 192, [ 16 ], [ 6 ],
 [ 2 ]
191, 192, [ 16 ], [ 2 ],
 [ 2 ]
192, 192, [ 8, 8 ], [ 2 ],
 []
193, 192, [ 16 ], [ 2, 2, 2 ],
 [ 2 ]
194, 320, [ 16 ], [ 4 ],
 []
195, 384, [ 16 ], [ 2, 2 ],
 [ 2 ]
196, 960, [ 16 ], [],
 []
197, 1920, [ 16 ], [ 2 ],
 []
\end{verbatim}}
\end{ttt}

\noindent
The list shows in every second line, in this order, the number of the conjugacy class in the numbering chosen by {\tt magma}, the order of the corresponding subgroups, the orbit type on the 16~lines the subgroups yield, and the type of the abelian quotient. The resulting group cohomology group is indicated on the following~line.

\begin{rems}
\begin{iii}
\item
To summarise, we find that
$H^1(H, D_5^*)$
is
\begin{abc}
\item[$\bullet$ ]
$0$
in 59 cases,
\item[$\bullet$ ]
$\bbZ/2\bbZ$
in 62 cases,
\item[$\bullet$ ]
$(\bbZ/2\bbZ)^2$
in 44 cases,
\item[$\bullet$ ]
$(\bbZ/2\bbZ)^3$
in 16 cases,
\item[$\bullet$ ]
$(\bbZ/2\bbZ)^4$
in three cases,
\item[$\bullet$ ]
$\bbZ/4\bbZ$
in nine cases,
\item[$\bullet$ ]
$\bbZ/2\bbZ \times \bbZ/4\bbZ$
in three cases, and
\item[$\bullet$ ]
$(\bbZ/2\bbZ)^2 \times \bbZ/4\bbZ$
in one case.
\end{abc}
In comparison with the Brauer groups of proper degree four del Pezzo surfaces, which may be only
$0$,
$\bbZ/2\bbZ$,
or
$(\bbZ/2\bbZ)^2$
\cite[Section~31, Table~3]{Ma} (see also~\cite{SD}), these figures are unexpectedly~rich.
\item
The running time was a little less than one second on one core of a AMD Phenom II X4 955 processor, using {\tt magma}, version 2.20.8.
\item
Note~that Manin's formula \cite[Proposition~31.3]{Ma} (see also \cite[Chapter~III, Proposition~8.18]{Ja}) does not apply to the computation of
$H^1(H, D_5^*)$,
as
$D_5^*$
is not a unimodular~lattice.
\end{iii}
\end{rems}

In the theorem below, we list some cases of particular interest.

\begin{theo}
\label{algBr}
Let\/~$k$
be a number field,
$\smash{X \subset \Pb^4_k}$
a degree four del Pezzo surface,
$H \subset \Pb^4_k$
be a\/
\mbox{$k$-rational}
hyperplane such that\/
$H \cap X$
is geometrically irreducible, and put\/
$U := X \!\setminus\! H$.

\begin{abc}
\item
Then\/
$\Br_1(U)/\Br_0(U)$
is isomorphic to\/
$0$,
$\bbZ/2\bbZ$,
$(\bbZ/2\bbZ)^2$,
$(\bbZ/2\bbZ)^3$,
$(\bbZ/2\bbZ)^4$,
$\bbZ/4\bbZ$,
$\bbZ/2\bbZ \!\times\! \bbZ/4\bbZ$,
or\/~$(\bbZ/2\bbZ)^2 \!\times\! \bbZ/4\bbZ$.
\item
Suppose~that
\begin{iii}
\item
the Galois group faithfully operating on the 16~lines
on\/~$X$
is isomorphic to the full\/
$W(D_5)$
or its index two subgroup.
Then~$\Br_1(U)/\Br_0(U) = 0$.
\item
the Galois group faithfully operating on the 16~lines
on\/~$X$
is isomorphic to the index five subgroup in\/
$W(D_5)$.
Then~$\Br_1(U)/\Br_0(U) = \bbZ/2\bbZ$.
\item
two of the five degenerate quadrics in the pencil associated with\/
$X$
are defined
over\/~$k$
and the Galois group faithfully operating on the 16~lines
on\/~$X$
is of
index\/~$20$
in\/
$W(D_5)$.
Then~$\Br_1(U)/\Br_0(U) = (\bbZ/2\bbZ)^2$.
\end{iii}
\end{abc}\medskip

\noindent
{\bf Proof.}
{\em
a), b.i) and b.ii) are immediate consequences from the list~above.\smallskip

\noindent
b.iii)
According to Facts~\ref{geomPic}.iii and~v), there is a canonical isomorphism
$$W(D_5) \longrightarrow (\bbZ/2\bbZ)^4 \rtimes S_5 \, ,$$
where the factor
$S_5$
is the group permuting the five degenerate quadrics in the pencil associated
with~$X$.
Thus,~by our first assumption, the Galois group
$H \hookrightarrow W(D_5)$
faithfully operating on the 16~lines is mapped isomorphically onto a subgroup
of~$(\bbZ/2\bbZ)^4 \rtimes S_3 \subset (\bbZ/2\bbZ)^4 \rtimes S_5$,
which is of
order~$96$.
On~the other hand,
order~$96$
is assumed, so that
$H$
goes over exactly into that semidirect~product. It~is not hard to check that 181 is the corresponding number in the~list.
}
\eop
\end{theo}

\subsubsection*{The transcendental part of the Brauer group.}

\begin{theo}
\label{trBr}
Let\/~$k$
be an algebraically closed field of
characteristic\/~$0$
and\/
$U := X \!\setminus\! H$,
for\/
$\smash{X \subset \Pb^4_k}$
a degree four del Pezzo surface and\/
$H \subset \Pb^4_k$
a hyperplane such that\/
$D := H \cap X$
is non-singular. Then there is a canonical~isomorphism
$$\Br(U) \cong H^1_\et(D,\bbQ/\bbZ) \, .$$
\noindent
{\bf Proof.}
{\em
Let
$n$
be any positive~integer. Then,~the Gysin sequence \cite[Exp.\ VII, Corollaire~1.5]{SGA5} (or~\cite[Chapter~VI, Remark~5.4.b)]{Mi}) reads
$$H^0_\et(D,\bbZ/n\bbZ) \longrightarrow H^2_\et(X,\mu_n) \longrightarrow H^2_\et(U,\mu_n) \longrightarrow H^1_\et(D,\bbZ/n\bbZ) \longrightarrow H^3_\et(X,\mu_n) \, .$$
Here,~as
$X$
is a del Pezzo surface, one has
$H^2_\et(X,\mu_n) \cong \Pic(X)/n\Pic(X)$
and
$H^3_\et(X,\mu_n) = 0$.
Moreover,
$H^0_\et(D,\bbZ/n\bbZ) = \bbZ/n\bbZ$
and the generator is mapped to the anticanonical class
in~$\Pic(X)/n\Pic(X)$.
Therefore,~in view of \cite[Proposition~II.6.5.c)]{Hart}, the sequence simplifies~to
$$0 \longrightarrow \Pic(U)/n\Pic(U) \longrightarrow H^2_\et(U,\mu_n) \longrightarrow H^1_\et(D,\bbZ/n\bbZ) \longrightarrow 0 \, ,$$
which shows that
$\Br(U)_n \cong H^1_\et(D,\bbZ/n\bbZ)$.
As~$n$
was arbitrary, this implies that
$\Br(U)_\tors \cong H^1_\et(D,\bbQ/\bbZ)$.

However,~for non-singular schemes in general, it is known that the cohomological Brauer group is torsion~\cite[Proposition~1.4]{GrBrII}.
}
\eop
\end{theo}

\begin{coro}
\label{trBrcor}
Let\/~$k$
be a number field and\/
$U := X \!\setminus\! H$,
for\/
$\smash{X \subset \Pb^4_k}$
a degree four del Pezzo surface and\/
$H \subset \Pb^4_k$
a\/
\mbox{$k$-rational}
hyperplane such that\/
$D := H \cap X$
is non-singular. Then,~for every
$n\geq 2$,
there is a canonical~monomorphism
\begin{equation}
\label{elltors}
[\Br(U)/\Br_1(U)]_n \hookrightarrow \left( J(D)(\bar k)_n\otimes_{\bbZ/n\bbZ} \mathbf{\mu}_n^{\vee} \right)^{\Gal(\overline{k}/k)} ,
\end{equation}
for\/
$J(D)$
the Jacobian variety
of\/~$D$.\medskip

\noindent
{\bf Proof.}
{\em
By~\ref{Brsubq}.iii), there is a canonical injection
$\Br(U)/\Br_1(U) \hookrightarrow \Br(U_{\overline{k}})^{\Gal(\overline{k}/k)}$.
Moreover,~Theorem~\ref{trBr} yields a canonical isomorphism
$$\Br(U_{\overline{k}})_n \cong H^1_\et(D_{\overline{k}}, \bbZ/n\bbZ) \cong H^1_\et(D_{\overline{k}}, \mu_n)\otimes_{\bbZ/n\bbZ} \mu_n^{\vee} \,.$$
Here,~$\smash{H^1_\et(D_{\overline{k}},\mathbf{\mu}_n) \cong \Pic(D_{\overline{k}})_n}$
by virtue of \cite[Chapter~III, Proposition~4.11]{Mi}. Moreover,~as
$D$
is a non-singular curve, one~has
$\smash{\Pic(D_{\overline{k}})_n \cong J(D)(\overline{k})_n}$,
according to the Picard interpretation of the~Jacobian.
}
\eop
\end{coro}

\begin{rems}
\begin{iii}
\item For the rest of the paper, we only consider
$2$-torsion
classes. In~this case, the inclusion (\ref{elltors}) simplifies to
$$[\Br(U)/\Br_1(U)]_2 \hookrightarrow  J(D)( k)_2\, .$$
Moreover, if $D$ has a \mbox{$k$-rational} point $O$ then we may choose this point to convert $D$ into an elliptic curve and get an injection into the $k$-rational $2$-torsion points of the elliptic curve $D$.  
\item
As~indicated for a more general situation in paragraph \ref{Brsubq}.iii), it may well be a hard problem to decide which elements on the right hand side of (\ref{elltors}) belong to the image of this monomorphism. In~some particular cases, however, we are able to write down transcendental Brauer classes, explicitly. Cf.~Theorem \ref{trBrexpl}~below.
\item
By Poincar\'e duality~\cite[Chapter~VI, Corollary~11.2]{Mi}, we have a perfect, Galois-invariant paring
$$H^1_\et(D_{\bar k}, \mu_n)\times H^1_\et (D_{\bar k}, \bbZ/n\bbZ) \rightarrow H^2_\et(D_{\bar k}, \mathbf{\mu}_n) \cong \bbZ/n\bbZ \,,$$
which induces an isomorphism
$$H^1_\et(D_{\bar k}, \tfrac{1}{n}\bbZ/\bbZ) \cong \Hom (H^1_\et(D_{\bar k}, \mu_n), \tfrac{1}{n}\bbZ/\bbZ) \cong \Hom(J(D)(\bar k)_n, \bbQ/\bbZ) \,.$$
Hence one obtains, similarly as in Corollary \ref{trBrcor}, a monomorphism
$$[\Br(U)/\Br_1(U)]_n \hookrightarrow \Hom_{\Gal(\overline{k}/k)} \!\left( J(D)(\bar k)_n, \bbQ/\bbZ\right) .$$
In~particular, a proper
\mbox{$n$-torsion}
class in
$\Br(U)/\Br_1(U)$
causes a
\mbox{$k$-rational}
\mbox{$n$-isogeny}
$J(D) \twoheadrightarrow D'$
to an elliptic curve~$D'$
with a
\mbox{$k$-rational}
proper
\mbox{$n$-torsion}~point.
\end{iii}
\end{rems}

\section{Explicit Brauer classes}
\label{sec_explBr}

\subsubsection*{The tame symbol.}\leavevmode\medskip

\noindent
Let~$X$
be an irreducible, non-singular scheme of dimension
$\leq \!2$ over a field $k$ of characteristic $0$.
Then~there is a canonical monomorphism
$\Br(X) \hookrightarrow \Br(k(X))$.
Moreover,~within
$\Br(k(X))$,
one has~\cite[Proposition~2.3]{GrBrII}
$$\Br(X) \cong  \!\bigcap_{x \in X^{(1)}}\!\! \Br(\Spec \calO_{X,x})\, ,$$
where
$X^{(1)}$
denotes the set of all codimension-one points on the
scheme~$X$.

In~other words, a Brauer class of the function field
$k(X)$
extends to the whole
of~$X$
if and only if it is unobstructed at any prime~divisor. And if this is the case then the extension is~unique.

The obstruction at a prime divisor
$x$
is tested by the so-called {\em ramification homomorphism}. This is the following composition~\cite[Proposition~2.1]{GrBrIII},
\begin{align*}
\ram_x \colon \Br(k(X)) \stackrel{\cong}{\longleftarrow} H^2_\et(X, \bbG_{m,k(X)}) \longrightarrow H^2_\et(X, i_* \bbZ_x) &\stackrel{\cong}{\longleftarrow} H^1_\et(\Spec k(x), \bbQ/\bbZ) \\
&\hspace{-16pt} = H^1(\Gal(\overline{k(x)}/k(x)), \bbQ/\bbZ) \, .
\end{align*}
Here,~the first map is an isomorphism due to~\cite[Chapter III, Corollary 5.6]{Ja}, the second is
the homomorphism induced by the valuation map
$$\val_x\colon \bbG_{m,k(X)} \rightarrow i_* \bbZ_x \, ,$$
and the third is the connecting homomorphism induced by the short exact sequence
$0 \rightarrow \bbZ \rightarrow \bbQ \rightarrow \bbQ/\bbZ \rightarrow 0$.\medskip

A way to explicitly write down a Brauer class of the function field is provided by the cyclic~algebras. For~these, let
$n \in \bbN$
and suppose that
$k(X)$
contains a primitive
\mbox{$n$-th}
root of
unity~$\zeta_n$.
Then,~for
$f,g \in k(X)^*$,
$$(f,g;\zeta_n) := \!\!\bigoplus_{0 \leq a,b < n}\!\!\! u^av^bk(X) \qquad {\rm for~} u^n = f, \,v^n = g, \,vu = uv\zeta_n \, ,$$
is an Azumaya algebra~\cite[\S15.4, Proposition]{Pi} and therefore defines a class in
$\Br(k(X))$.
The~corresponding cohomology class
in~$\smash{H^2(\Gal(\overline{k(X)}/k(X)), \overline{k(X)}{}^*)}$
may be described by the~cocycle
$\Phi$
such that~\cite[\S15.1, Proposition~a]{Pi}
$$\Phi(\sigma, \tau) = \left\{
\begin{array}{lc}
1 & \;\;{\rm if}\; i+j < n \, ,\\
f & \;\;{\rm if}\; i+j \geq n \, ,
\end{array}
\right.$$
for
$0 \leq i,j < n$
chosen by the condition that
$\sigma(\!\sqrt[n]{\mathstrut g}) = \zeta_n^i \sqrt[n]{\mathstrut g}$
and~$\tau(\!\sqrt[n]{\mathstrut g}) = \zeta_n^j \sqrt[n]{\mathstrut g}$.\medskip

Finally,~recall that there is the tame symbol
$$\textstyle \delta_x\colon k(X)^* \times k(X)^* \longrightarrow k(x)^*, \quad (f,g) \mapsto (-1)^{\nu_x(f)\nu_x(g)} \overline{(\frac{f^{\nu_x(g)}}{g^{\nu_x(f)}})} \, ,$$
where the bar denotes the natural residue map
$\calO_{X,x} \to k(x)$.

\begin{prop}
Let\/~$n \in \bbN$
and fix a primitive\/
\mbox{$n$-th}
root of unity\/
$\zeta_n \in k(X)^*$.
Then~the following diagram commutes,
$$
\diagram
k(X)^* \!\times k(X)^* \rrto^{\;\;\;\;\;\;\;\;\;\;\delta_x} \dto_{R_{n,\zeta_n}} && k(x)^* \dto^\pi \\
\Br(k(X)) \rrto^{\!\!\!\!\!\!\!\!\!\!\!\!\!\!\!\!\!\!\!\!\!\!\!\!\ram_x} && H^1(\Gal(\overline{k(x)}/k(x)), \bbQ/\bbZ) \, .
\enddiagram
$$
Here, the downward arrow\/
$R_{n,\zeta_n}$
to the left sends\/
$(f,g)$
to the Brauer class of\/
$(f,g;\zeta_n)$,
while the one to the right is the canonical map
$$k(x)^* \twoheadrightarrow k(x)^*/(k(x)^*)^n \stackrel{\cong}{\longrightarrow} H^1(\Gal(\overline{k(x)}/k(x)), \mu_n) \, ,$$
followed by the isomorphism induced by the identification\/
$\mu_n \rightarrow \frac1n\bbZ/\bbZ$,
$\zeta_n \mapsto -\frac1n$.\medskip

\noindent
{\bf Proof.}
{\em If
$\nu_x(f) = \nu_x(g) = 0$
then
$\delta_x(f,g) = 1$
such that
$(f,g)$
is mapped to
$0$
via the upper right~corner.

On~the other hand, the cohomology class of
$R_{n,\zeta_n}(f,g) = (f,g;\zeta_n)$
comes via inflation from
$\smash{H^2(\Gal(k(X)(\!\sqrt[n]{\mathstrut g})/k(X)), \overline{k(X)}{}^*)}$
and is represented by a cocycle whose image
is~$\{1,f\}$.
Since~the field extension
$\smash{k(X)(\!\sqrt[n]{\mathstrut g})/k(X)}$
is unramified
at~$x$
and one has
$\nu_x(1) = \nu_x(f) = 0$,
it follows that
$\ram_x(f,g;\zeta_n)$
is represented by the cocycle that is constantly~zero.
Hence,~$\ram_x(f,g;\zeta_n) = 0$.
Thus,~via the lower left~corner,
$(f,g)$
is mapped
to~$0$,
too.

Next,~we observe that the maps
$\delta_x$
and
$R_{n,\zeta_n}$
are bimultiplicative and antisymmetric. Thus,~in order to complete the proof, it suffices to consider the case that
$\nu_x(f) = 1$
and
$\nu_x(g) = 0$.\smallskip

\noindent
{\em Upper right path:}
We have
$\delta_x(f,g) = \overline{(1/g)}$,
which is mapped
under~$\pi$
to the class of the homomorphism
$$\sigma \mapsto \frac{\sigma(\!\sqrt[n]{\overline{(1/g)}})}{\sqrt[n]{\overline{(1/g)})}} = \frac{\sqrt[n]{\overline{g}}}{\sigma(\!\sqrt[n]{\overline{g}})}$$
in~$H^1(\Gal(k(x)(\!\sqrt[n]{\overline{g}})/k(x)), \mu_n) \hookrightarrow H^1(\Gal(\overline{k(x)}/k(x)), \mu_n)$.
In~other words,
$\pi(\delta_x(f,g))$
maps the element
$\sigma \in \Gal(k(x)(\!\sqrt[n]{\overline{g}})/k(x))$
such that
$\sqrt[n]{\overline{g}} \mapsto \zeta_n^i \sqrt[n]{\overline{g}}$
to
$\zeta_n^{-i}$,
which finally gets identified
with~$\frac{i}n$.\smallskip

\noindent
{\em Lower left path:}
Here,~$(f,g)$
is sent first
by~$R_{n,\zeta_n}$
to
$(f,g;\zeta_n)$,
which is described by the~cocycle
$\Phi$
such that
$$\Phi(\sigma, \tau) = \left\{
\begin{array}{lc}
1 & \;\;{\rm if}\; i+j < n \, ,\\
f & \;\;{\rm if}\; i+j \geq n \, ,
\end{array}
\right.$$
for
$0 \leq i,j < n$
chosen in the way that
$\sigma(\!\sqrt[n]{\mathstrut g}) = \zeta_n^i \sqrt[n]{\mathstrut g}$
and~$\tau(\!\sqrt[n]{\mathstrut g}) = \zeta_n^j \sqrt[n]{\mathstrut g}$.
Applying~the
valuation~$\nu_x$,
one finds
$$(\nu_x \!\circ\! \Phi)(\sigma, \tau) = \left\{
\begin{array}{lc}
0 & \;\;{\rm if}\; i+j < n \, ,\\
1 & \;\;{\rm if}\; i+j \geq n \, .
\end{array}
\right.$$
This,~however, is exactly the image~of
$$\textstyle \Psi\colon \sigma \mapsto \frac{i}n \, , \qquad {\rm for~} \sigma(\!\sqrt[n]{\mathstrut g}) = \zeta_n^i \sqrt[n]{\mathstrut g} \, ,$$
understood as an element
in~$\smash{H^1(\Gal(k(X)(\sqrt[n]{\mathstrut g})/k(X)), \bbQ/\bbZ)}$,
under the connecting homomorphism associated with
$0 \rightarrow \bbZ \rightarrow \bbQ \rightarrow \bbQ/\bbZ \rightarrow 0$.
Thus,~$\smash{\ram_x(R_{n,\zeta_n}(f,g))}$
maps the element
$\smash{\sigma \in \Gal(k(x)(\!\sqrt[n]{\overline{g}})/k(x))}$
sending
$\smash{\sqrt[n]{\overline{g}} \mapsto \zeta_n^i \sqrt[n]{\overline{g}}}$
to~$\smash{\frac{i}n}$.
This~completes the~proof.
}
\eop
\end{prop}

\begin{coro}
Let\/~$X$
be an irreducible, non-singular scheme of dimension\/
$\leq \!2$ over a field $k$ of characteristic $0$,
$n \in \bbN$,
$\zeta_n \in k(X)$
a primitive
\mbox{$n$-th}
root of unity,
$f,g \in k(X)^*$,
and\/
$x \in X^{(1)}$
a prime~divisor. Then~the Brauer class
of\/~$(f,g;\zeta_n)$
is unobstructed
at\/~$x$
if and only if\/
$\delta_x(f,g)$
is an
\mbox{$n$-th}
power
in\/~$k(x)$.
\eop
\end{coro}

\begin{rems}
\begin{iii}
\item
In~Theorem~\ref{trBr}, the Gysin homomorphism appeared instead of the tame~symbol. It~is, in fact, true that these two maps are compatible. Cf.,~for example, \cite[Lemma~2.5]{BMT}.
\item
In \cite[\S3]{AM}, M.~Artin and D.~Mumford provide a more complete picture of the ramification~homomorphism. For example,~the following is~shown.

Let~$X$
be a proper surface over an algebraically closed field and
$U := X \!\setminus\! D$
for
$D$
an irreducible non-singular~curve. Suppose~that
$(f,g;\zeta_n)$
defines a Brauer class
on~$U$.
Then~the divisor
$\div t$
of the rational~function
$t := \delta_x(f,g) \in k(D)$,
for
$x$
the generic point
of~$D$,
is divisible
by~$n$.

Of~course, on a curve of genus
$\geq \!1$,
this is not enough to imply that
$t$
has to be an
\mbox{$n$-th}~power.
\end{iii}
\end{rems}

\subsubsection*{Two model cases.}

\begin{theo}[An algebraic Brauer class]
\label{algBrexpl}
Let\/~$k$
be any field of characteristic $0$,
$\smash{X \subset \Pb^4_k}$
a del Pezzo surface of degree four
over\/~$k$,
and\/
$U := X \!\setminus\! H$
for\/
$H := V(X_4) \subset \Pb^4_k$.
Suppose~that the pencil of quadrics associated with\/
$X$
contains a\/
\mbox{$k$-rational}
quadric of
rank\/~$4$
that may be written in the~form
$$l_1l_2 - l_3^2 + dl_4^2 \, ,$$
for linear forms\/
$l_1, \ldots, l_4$
and\/
$d \in k^*$
a non-square.

\begin{iii}
\item
Then
$$\textstyle (\frac{l_1}{X_4}, d; -1)$$
defines an algebraic Brauer~class\/
$\alpha \in \Br_1(U)_2$.
\item
If\/~$V(X_4) \cap X$
is reduced and has a\/
\mbox{$k$-rational}
point then\/
$\alpha$
does not extend
to\/~$X$
and is, in particular, nontrivial.
\end{iii}\medskip

\noindent
\mbox{{\bf Proof.} {\em i)}}
{\em
Let
$D \subset U$
be any prime divisor and
$x \in U^{(1)}$
its generic~point.
Then~$\nu_x(d) = 0$
and
$\smash{\nu_x(\frac{l_1}{X_4}) = 0}$,
unless
$D$
is a component
of~$V(l_1) \cap X$.
Thus,~if
$D \not\subset V(l_1) \cap X\!$
then
$\!\smash{\delta_x(\frac{l_1}{X_4}, d) = 1}$.

On~the other hand, if
$D \subset V(l_1) \cap X$
then
$\smash{\delta_x(\frac{l_1}{X_4}, d) = \frac1{d^{\nu_x(l_1)}}}$,
which is a square on
$D$,
since the equation
$- l_3^2 + dl_4^2 = 0$
implies
$$\textstyle d = \frac{l_3^2}{l_4^2} \, .$$
The~Brauer class
$\alpha$
is clearly algebraic, as it gets annihilated under base extension
to~$k(\sqrt{d})$.\smallskip

\noindent
ii)
Take~now
$D$
to be a component of
$V(X_4) \cap X$.
Then~$\smash{\delta_x(\frac{l_1}{X_4}, d) = d}$
for
$x \in U^{(1)}$
the generic point
of~$D$,
which cannot be a square when 
$D$
has a
\mbox{$k$-rational}
point.%
}%
\eop
\end{theo}

\begin{theo}[A transcendental Brauer class]
\label{trBrexpl}
Let\/~$k$
be any field of characteristic $0$ and\/
$\smash{X \subset \Pb^4_k}$
a del Pezzo surface of degree four
over\/~$k$
that is given by a system of equations of the~type
\begin{align*}
l_1l_2 + au^2 &= X_4l_3 \, ,\\
l_3l_4 + bv^2 &= X_4l_1 \, ,
\end{align*}
for linear forms\/
$l_1, \ldots, l_4, u, v$,
and\/
$a,b \in k^*$.
Assume~that the forms\/
$l_1$,
$l_3$,
$u$,
and\/~$v$
are linearly~independent.\smallskip

\noindent
Put\/
\mbox{$U := X \!\setminus\! H$}
for\/
$H := V(X_4) \subset \Pb^4_k$.

\begin{iii}
\item
Then
$$\textstyle (\frac{bl_1}{X_4}, \frac{al_3}{X_4}; -1)$$
defines a Brauer~class\/
$\tau \in \Br(U)_2$.
\item
If\/~$D := V(X_4) \cap X$
is geometrically integral and\/
$\smash{\frac{l_1}{l_3}}$
is not the square of a rational function
on\/~$\smash{D_{\overline{k}}}$
then\/
$\smash{\tau_{U_{\overline{k}}}}$
does not extend to\/
$\smash{X_{\overline{k}}}$.
In~particular,
$\tau$
is~transcendental.
\end{iii}\medskip

\noindent
{\bf Proof.}
{\em
i)
Let
$D \subset U$
be any prime divisor and
$x \in U^{(1)}$
its generic~point. Then
$\smash{\nu_x(\frac{bl_1}{X_4}) = 0}$
and
$\smash{\nu_x(\frac{al_3}{X_4}) = 0}$,
unless
$D \subset V(l_1) \cap X$
or~$D \subset V(l_3) \cap X$.
Thus,~if
$D$
is not a component of either of these two subschemes then
$\smash{\delta_x(\frac{bl_1}{X_4}, \frac{al_3}{X_4}) = 1}$.

On~the other hand, if
$D \subset V(l_1) \cap X$,
but
$D \not\subset V(l_3) \cap X$,
then one finds that
$\smash{\delta_x(\frac{bl_1}{X_4}, \frac{al_3}{X_4}) = (\frac{X_4}{al_3})^{\nu_x(l_1)}}$,
which is a square
on~$D$,
since the equation
$au^2 = X_4l_3$
implies
$$\textstyle \frac{X_4}{al_3} = \frac{u^2}{l_3^2} \, .$$
Similarly,~if
$D \subset V(l_3) \cap X$,
but
$D \not\subset V(l_1) \cap X$,
then
$\smash{\delta_x(\frac{bl_1}{X_4}, \frac{al_3}{X_4}) = (\frac{bl_1}{X_4})^{\nu_x(l_3)}}$,
which is a square
on~$D$,
since
$bv^2 = X_4l_1$
yields
$$\textstyle \frac{bl_1}{X_4} = \frac{l_1^2}{v^2} \, .$$
Finally,~there cannot be a divisor that is contained in both,
$V(l_1) \cap X$
and
$V(l_3) \cap X$,
as this would imply
$l_1 = l_3 = u = v = 0$,
which, according to our assumptions, is fulfilled only by a single~point.\smallskip

\noindent
ii)
Let~now
$x$
be the generic point of
$V(X_4) \cap X$.
Then~$\smash{\delta_x(\frac{bl_1}{X_4}, \frac{al_3}{X_4}) = \frac{al_3}{bl_1}}$.
Since
$\smash{\frac{a}b}$
is a square
in~$\smash{\overline{k}}$,
this immediately implies the~claim.
}
\eop
\end{theo}

\subsubsection*{Constancy of the local evaluation maps.}

\begin{rem}
Let~now
$k$
be a number field and
$U$
be a
\mbox{$k$-scheme},
on which a Brauer class
$\alpha \in \Br(U)$
is given by
$(f,g;-1)$,
for rational functions
$f,g \in k(U)$.
Moreover,~let
$\nu$
be any place
of~$k$
and
$x \in U(\bbQ_\nu)$,
$x \not\in \supp\div f \cup \supp\div g$,
a~point. Then
$$\ev_{\alpha,\nu}(x) = \left\{
\begin{array}{ll}
      0 & \;\;{\rm if}\; (f(x),g(x))_\nu = 1 \, ,\\
\frac12 & \;\;{\rm if}\; (f(x),g(x))_\nu = -1 \, ,
\end{array}
\right.
$$
where
$(.,.)_\nu \colon k_\nu^* \times k_\nu^* \rightarrow \{\,1,-1\,\}$
denotes the usual
\mbox{$\nu$-adic}
Hilbert~symbol.

Indeed,~according to the definition,
$\ev_{\alpha,\nu}(x) = x^*\alpha$
and our assumptions guarantee that the pull-back may be obtained in a naive way as the cyclic (quaternion) algebra
$(f(x),g(x);-1)$
over~$k_\nu$.
Finally,~the invariant of such an algebra is directly related, by the above rule, to the Hilbert symbol~\cite[\S18.4, Exercise~4]{Pi}.
\end{rem}

\begin{lem}
\label{algconst}
Let\/~$k$
be a number field,
$\smash{\calX \subset \Pb^4_{\calO_k}}$
the closed subscheme defined by two quadratic forms with coefficients
in\/~$\calO_k$,
and\/
$\calU := \calX \!\setminus\! \calH$
for\/
$\calH := V(X_4) \subset \Pb^4_{\calO_k}$.
Suppose~that the generic fibre\/
$X$
of\/
$\calX$
is a non-singular surface and let\/
$U$
and\/~$\alpha \in \Br(U)_2$
be as in Theorem~\ref{algBrexpl}.\smallskip

\noindent
Assume that one of the quadratic forms defining\/
$\calX$
may be written in the~form
$$l_1l_2 - l_3^2 + dl_4^2 \, ,$$
for\/
$l_1, \ldots, l_4$
linear forms having coefficients
in\/~$\calO_k$
and\/
$d \in \calO_k$
a non-square.
Furthermore,~let\/
$\frakp$
be a prime ideal in\/
$\calO_k$
such
that\/~$\frakp \notd 2d$,

\begin{iii}
\item[$\bullet$ ]
the linear forms
$(l_1 \bmod \frakp), \ldots, (l_4 \bmod \frakp)$
are\/
\mbox{$\calO_k/\frakp$-linearly}
independent, and
\item[$\bullet$ ]
the reduction modulo\/
$\frakp$
of the cusp given by\/
$l_1 = \ldots = l_4 = 0$
does not lie on the
reduction~$\calU_\frakp$.
\end{iii}
Then~the local evaluation~map
$$\ev_{\alpha,\nu_\frakp}\colon \calU(\calO_{k,\frakp}) \longrightarrow \bbQ/\bbZ$$
is constantly~zero.\medskip

\noindent
{\bf Proof.}
{\em
Let~us denote integral points
$x \in \calU(\calO_{k,\frakp})$
in the form
$x = (x_0:\ldots:x_3:1)$,
for
$x_0, \ldots, x_3 \in \calO_{k,\frakp}$.
The local evaluation~map
$\ev_{\alpha,\nu_\frakp}$
is known to be locally constant with respect to the
\mbox{$\frakp$-adic}
topology
on
$U(k_{\nu_\frakp})$~\cite[Chapter~IV, Proposition~2.3.a.ii)]{Ja}.
Consequently,~it suffices to prove the assertion for points
$x \in \calU(\calO_{k,\frakp})$
such that
$l_1(x) \neq 0$
and~$l_2(x) \neq 0$.
Thus,~what we have to show~is
$$(l_1(x_0, \ldots, x_3, 1), d)_{\nu_\frakp} = 1 \, ,$$
for every point
$x = (x_0:\ldots:x_3:1) \in U(\calO_{k,\frakp})$
satisfying the inequalities~above.

For~this, we first observe that the equation given yields
$l_1(x)l_2(x) = l_3^2(x) - dl_4^2(x)$,
which is a norm
from~$\smash{k_{\nu_\frakp}(\sqrt{d})}$.
Therefore,~in the Hilbert symbol, we may replace
$l_1(x)$
by~$l_2(x)$.

If~$(d \bmod \frakp) \in \calO_k/\frakp$
is a square then assertion is certainly~true. Thus~assume that
$(d \bmod \frakp) \in \calO_k/\frakp$
is a non-square.
Then~$(l_1(x), d)_{\nu_\frakp} = -1$
would mean that
$\nu_\frakp(l_1(x))$
is~odd. In~which case,
$\nu_\frakp(l_2(x))$
is also~odd. In~particular,
$\frakp | l_1(x)$
and~$\frakp | l_2(x)$,
which directly implies that
$\frakp | l_3(x)$
and
$\frakp | l_4(x)$,~too.

However,~according to our assumptions, such an
\mbox{$\calO_{k,\frakp}$-valued}
point
on~$\calU$
does not~exist. Indeed,~its reduction
modulo~$\frakp$
would coincide with that of the cusp, which is supposed not to lie
on~$\calU_\frakp$.
}
\eop
\end{lem}

\begin{rem}
In~particular, we see that if
$\frakp \notd 2$
is a prime ideal of good reduction
of~$\calX$
then
$\ev_{\alpha,\frakp}$
is constantly~zero.

Indeed,~then
$\calX_\frakp$
is a degree four del Pezzo surface
over~$\bbF_{\!\frakp}$.
Hence,~the rank of a degenerate quadric in the pencil cannot drop
below~$4$~\cite[Proposition 3.26.iv)]{Wi}.
This~implies immediately that
$\frakp \notd d$,
as well as linear independence of
\mbox{$(l_1 \bmod \frakp), \ldots, (l_4 \bmod \frakp)$}.
Moreover,~the reduction of the cusp would be a singular point
on~$\calU_\frakp$,
which, by assumption, does not~exist.
\end{rem}

\begin{lem}
\label{trconst}
Let\/~$k$
be a number field and\/
$\smash{\calX \subset \Pb^4_{\calO_k}}$
be the closed subscheme defined by a system of equations of the~type
\begin{align*}
l_1l_2 + au^2 &= X_4l_3 \, ,\\
l_3l_4 + bv^2 &= X_4l_1 \, ,
\end{align*}
for linear forms\/
$l_1, \ldots, l_4, u, v$
having coefficients
in\/~$\calO_k$
and\/
$a,b \in \calO_k^*$.
Assume~that the forms\/
$l_1$,
$l_3$,
$u$,
and\/~$v$
are
\mbox{$k$-linearly}
independent and that the generic fibre\/
$X$
of\/~$\calX$
is a non-singular~surface.\smallskip

\noindent
Put\/
\mbox{$\calU := \calX \!\setminus\! \calH$}
for\/
$\calH := V(X_4) \subset \Pb^4_{\calO_k}$
and
let\/~$U$
and\/~$\tau \in \Br(U)_2$
be as in Theorem~\ref{trBrexpl}.
Then,~for every prime ideal\/
$\frakp$
in\/~$\calO_k$
such that\/
$\frakp \notd 2ab$,
the local evaluation~map
$$\ev_{\tau,\nu_\frakp}\colon \calU(\calO_{k,\frakp}) \longrightarrow \bbQ/\bbZ$$
is constantly~zero.\medskip

\noindent
{\bf Proof.}
{\em
Let~us denote integral points
$x \in \calU(\calO_{k,\frakp})$
in the form
$x = (x_0:\ldots:x_3:1)$
for
$x_0, \ldots, x_3 \in \calO_{k,\frakp}$.
As~the local evaluation~map
$\ev_{\tau,\nu_\frakp}$
is locally constant with respect to the
\mbox{$\frakp$-adic}
topology
on
$U(k_{\nu_\frakp})$,
it suffices to prove the assertion for points
$x \in \calU(\calO_{k,\frakp})$
such that
$l_i(x) \neq 0$,
for~$i = 1,\ldots,4$.
Thus,~what we have to show~is
$$(bl_1(x), al_3(x))_{\nu_\frakp} = 1 \, ,$$
for every point
$x = (x_0:\ldots:x_3:1) \in \calU(\calO_{k,\frakp})$
satisfying the inequalities~above.

First~of all, observe that
$a$
and
$b$
are
\mbox{$\frakp$-adic}
units. Furthermore,~the first equation yields
$(-a)l_1(x)l_2(x) = a^2u^2(x) - al_3(x)$,
which is a norm
from~$\smash{k_{\nu_\frakp}(\!\sqrt{al_3(x)})}$.
Therefore,~in the Hilbert symbol, we may replace
$bl_1(x)$
by~$(-ab)l_2(x)$.

Analogously,~we have
$\smash{(-b)l_3(x)l_4(x) = b^2v^2(x) - bl_1(x)}$
being a norm from
$\smash{k_{\nu_\frakp}(\!\sqrt{bl_1(x)})}$.
Thus,~we may as well replace
$al_3(x)$
by~$(-ab)l_4(x)$.

This~shows, if one of the numbers
$l_1(x), l_2(x)$
is of even
\mbox{$\frakp$-adic}
valuation, as well as one of the numbers
$l_3(x), l_4(x)$,
then assertion is certainly~true.
Thus,~suppose without restriction that
$\nu_\frakp(l_1(x))$
and~$\nu_\frakp(l_2(x))$
are both~odd. We~have to show in this situation that
$al_3(x) \in k_{\nu_\frakp}^*$
is a~square.

For~this, as
$\nu_\frakp(l_1(x))$
is odd and
$\nu_\frakp(bv^2(x))$
is even, the second equation~implies
\begin{align*}
\nu_\frakp(l_3(x)l_4(x)) &= \min(\nu_\frakp(bv^2(x)), \nu_\frakp(l_1(x))) \\
                         &\leq \nu_\frakp(l_1(x)) \, ,
\end{align*}
hence
$$\nu_\frakp(l_3(x)) \leq \nu_\frakp(l_3(x)l_4(x)) \leq \nu_\frakp(l_1(x)) < \nu_\frakp(l_1(x)l_2(x)) \, .$$
Note~here that
$\nu_\frakp(l_2(x)) \geq 1$,
the valuation of
$l_2(x)$
being~odd.
Consequently,~the first equation implies that
$al_3(x)$
is a square
in~$k_{\nu_\frakp}^*$.
This~completes the~proof.
}
\eop
\end{lem}

\section{Examples}
\label{sec_ex1}

Recall~\cite[Definition~2.1]{DW} (cf.~\cite[\S7.1]{PR}) that {\em strong approximation off a certain set\/
$S$
of places\/} is said to hold for an algebraic variety
$U$
defined over a number
field~$k$
if the image of the set
$U(k)$
of
\mbox{$k$-rational}
points
on~$U$
is dense in the space
$U(\bbA_k^S)$
of adelic points
on~$U$
outside~$S$.

\subsubsection*{Del Pezzo surfaces of degree four.}\leavevmode\\[1mm]
In~our first example, only a transcendental Brauer class~occurs.

\begin{ex}
\label{ex1}
Let~$X \subset \Pb^4_\bbQ$
be given by the system of~equations
\begin{align*}
           X_0X_1+X_2^2 &= X_4X_3 \, , \\
X_3(2X_1+X_2+X_3)+X_0^2 &= X_4X_1
\end{align*}
and put
$U := X \!\setminus\! H$,
for the hyperplane
$\smash{H := V(X_4) \subset \Pb^4_\bbQ}$.

\begin{iii}
\item
Then~$X$
is a del Pezzo surface of degree~four.
\item
$D := H \cap X$
is a non-singular curve of genus~one. Furthermore,
$$D(\bbQ) = \{(0\!:\!1\!:\!0\!:\!0\!:\!0), (0\!:\!-1\!:\!0\!:\!2\!:\!0)\} \, .$$
\item
The manifold\/
$X(\bbR)$
is connected, its submanifold
$U(\bbR)$
is connected, too, and
$U$
is strongly unobstructed
at~$\infty$.
\item
One has
$\Br(U)/\Br_0(U) \cong \bbZ/2\bbZ$,
a representative of the nontrivial element being
$$\textstyle \tau\colon (\frac{X_1}{X_4}, \frac{X_3}{X_4}; -1) \, .$$
Moreover,~$\Br_1(U)/\Br_0(U) = 0$.
I.e.,~the Brauer class
$\tau$
is~transcendental.
\end{iii}\medskip

\noindent
{\bf Proof.}
i) and~ii)
{\tt magma} reports that both
$X$
and~$D$
are non-singular. More~precisely,
$D$~is
isomorphic to the elliptic curve
$E\colon y^2z + xyz = x^3 + 4xz^2$,
the Mordell-Weil group of which is isomorphic
to~$\bbZ/2\bbZ$.\smallskip

\noindent
iii)
The~pencil of quadrics
in~$\Pb^4$
corresponding
to~$X$
contains only one degenerate quadric that is defined
over~$\bbR$.
According to~\cite[Theorem~3.4 and Lemma~5.1]{VAV}, this implies
$\Br(X_\bbR)/\Br(\bbR) = 0$,
which is enough to show that
$X(\bbR)$
is connected \mbox{\cite[Th\'eor\`eme 1.4]{Si}}.

As~$X$
has no real point such that
$x_1 = x_4 = 0$,
the rational map
$\pi\colon X \ratarrow \Pb^1$,
given by
$x \mapsto x_1/x_4$,
is well-defined on real points.
Furthermore,
$U(\bbR) = \pi^{-1}(\Ab^1(\bbR))$.
A~calculation using Gr\"obner bases shows that
$\pi$
has exactly 10 singular fibres, among which only two, those over
$0$
and~$\infty$,
are~real. Finally,~experiments with
$t = 1$
and
$t = -1$
indicate that the fibres
$\pi^{-1}(t)$,
for
$t>0$
as well as for
$t<0$,
are genus-one curves with a connected set of real~points. As~connectedness
of~$\pi^{-1}(0)$
is directly checked, this is enough to guarantee that
$U(\bbR)$
is~connected.

Since~$U(\bbR)$
is connected and
$D$
is a non-singular, geometrically irreducible scheme having real points, Theorem~\ref{unobstr} shows that
$U$~is
strongly unobstructed at the real~place.\smallskip

\noindent
iv)
A~calculation using Gr\"obner bases shows that the Galois group operating on the 16~lines
on~$X$
is the full
$W(D_5)$
of
order~$1920$.
Therefore,~Theorem \ref{algBr}.b.i) immediately implies that
$\Br_1(U)/\Br_0(U) = 0$.

On~the other hand, we have
$\#D(\bbQ) = 2$.
Moreover,~$E$
has two
\mbox{$\bbQ$-rational}
\mbox{$4$-isogenies,}
but no nontrivial
\mbox{$n$-isogenies}
for any
other~$n \geq 3$.
As~the isogenous curves have no proper
\mbox{$\bbQ$-rational}
\mbox{$4$-torsion}
points, Corollary~\ref{trBrcor} shows that there is an injection
$\Br(U)/\Br_1(U) \hookrightarrow \bbZ/2\bbZ$.
Consequently,~$\Br(U)/\Br_0(U)$
is of order at
most~$2$.
As~far as a concrete description is asked for, Theorem~\ref{trBrexpl} applies. It~shows that
$\smash{(\frac{X_1}{X_4}, \frac{X_3}{X_4}; -1)}$
defines a Brauer class
on~$U$.

In~order to prove that this is indeed a nontrivial class, the simplest argument is probably that its local evaluation map 
$\ev_{\tau, \nu_2}$
is non-constant, which we show in Example~\ref{ex1b}.i)~below.
\eop
\end{ex}

\begin{ex}[continued]
\label{ex1b}
Let~$X$,
$U$,
and~$\tau$
be as~above. Moreover,~let
$\calX \subset \Pb^4_\bbZ$
be defined by the same system of equations
as~$X$
and put
$\calU := \calX \!\setminus\! \calH$,
for
$\calH := V(X_4) \subset \Pb^4_\bbZ$.
Then

\begin{iii}
\item
For all primes
$p \neq 2$
including the archimedean one, the local evaluation map
$\ev_{\tau, \nu_p}\colon \calU(\bbZ_p) \to \bbQ/\bbZ$
is constantly~zero. On~the other hand,
$\ev_{\tau, \nu_2}\colon \calU(\bbZ_2) \to \bbQ/\bbZ$
is non-constant.
\item
For~$p = 17$,
even the local evaluation map
$\ev_{\tau, \nu_{17}}\colon U(\bbQ_{17}) \to \bbQ/\bbZ$
is constantly zero.
\end{iii}\medskip

\noindent
{\bf Proof.}
i)
For~all primes with the exception of the infinite one, this follows directly from Lemma~\ref{trconst}. To~prove constancy for
$\ev_{\tau, \nu_\infty}$,
we have to show that, for a real point
$x = (x_0\!:\!\ldots\!:\!x_3\!:\!1) \in U(\bbR)$,
the coordinates
$x_1$
and
$x_3$
cannot both be~negative.

Assume~the contrary. In~that case, the first equation shows
$x_0 > 0$.
Moreover,
$x_2^2 + |x_3| = |x_0x_1|$
and hence
$x_2^2 \leq |x_0x_1|$,
implying
$x_2 \leq x_0$
or~$x_2 \leq -x_1$.
In~the first case,
\begin{align*}
x_3(2x_1+x_2+x_3)+x_0^2 &= 2x_1x_3+x_2x_3+x_3^2+x_0^2 \\
&\geq 2x_1x_3+x_0x_3+x_3^2+x_0^2 > 2x_1x_3+(x_0+x_3)^2 > 0 \, ,
\end{align*}
in contradiction to the second~equation. And~in the second case, the same is true, since
$x_2x_3 \geq -x_1x_3$.

Finally,~non-constancy of
$\ev_{\tau, \nu_2}$
may easily be seen directly. Indeed,~for the integral point
$x = (-1\!:\!-1\!:\!-1\!:\!2\!:\!1) \in \calU(\bbZ_2)$,
one has
$(x_1, x_3)_2 = (-1, 2)_2 = 1$,
hence
$\ev_{\tau, \nu_2}(x) = 0$.
On~the other hand, the
\mbox{$\bbZ_2$-valued}
point
$\smash{x' = (\frac25\!:\!\frac25\!:\!\frac15\!:\!\frac15\!:\!1) \in \calU(\bbZ_2)}$
yields
$(x_1, x_3)_2 = (\frac25, \frac15)_2 = -1$
and
$\ev_{\tau, \nu_2}(x) = \frac12$.\smallskip

\noindent
ii)
The local evaluation map
$\ev_{\tau, \nu_{17}}$
is locally constant with respect to the
\mbox{$17$-adic}
topology on
$U(\bbQ_{17})$.
Thus,~it suffices to prove the claim for points
$x = (x_0\!:\!\ldots\!:\!x_4) \in U(\bbQ_{17})$
such that
$x_1 \neq 0$
and~$x_3 \neq 0$.
Moreover,~the assumption implies
that~$x_4 \neq 0$.
Let~us assume that the coordinates are normalised such that
$x_0, \ldots, x_4 \in \bbZ_{17}$
and at least one is a~unit.

If~$x_4$
is a unit then
$x$
defines a
\mbox{$\bbZ_{17}$-valued}
point
on~$\calU$.
In~this case,
$\ev_{\tau, \nu_{17}}(x) = 0$,
as follows from~i).
Thus,~assume that
$17|x_4$.
Then~$x$
reduces modulo
$17$
to a point on the reduction
$D_{17}$
of the
curve~$D$.
A~direct inspection shows that
$D_{17}$
has exactly~17
$\bbF_{\!17}$-rational
points, that
$\xi_1\xi_3 \neq 0$
for exactly 14 of these, and that
$\xi_1/\xi_3 \in \bbF_{\!17}^*$
is a square for each of~them. Therefore,~if
$17 \notd x_1$
and
$17 \notd x_3$
then
$x_1/x_3 \in \bbQ_{17}^*$
is a~square and, consequently,
$\smash{(\frac{x_1}{x_4}, \frac{x_3}{x_4})_{17} = (-\frac{x_1}{x_3}, \frac{x_3}{x_4})_{17} = 1}$.

Moreover,~if
$17|x_4$
and
$17|x_3$
then the equations of the surface imply
$17|x_0$
and
$17|x_2$,
so that
$x_1$
is a~unit. As~the second equation shows that
$-x_3(2x_1+x_2+x_3)$
is a norm from
$\smash{\bbQ_{17}({\!\sqrt{x_1/x_4}})}$,
one has
$\smash{(\frac{x_1}{x_4}, \frac{x_3}{x_4})_{17} = (\frac{x_1}{x_4}, \frac{-2x_1-x_2-x_3}{x_4})_{17} = (\frac{x_1}{x_4}, \frac{-2x_1}{x_4})_{17} = 1}$.\vskip0.1mm

Finally, assume that
$17|x_4$
and
$17|x_1$.
Then the equations of the surface imply
$17|x_2$
and that
$x_0$
and
$x_3$
are units such that
$x_0^2 + x_3^2 \equiv 0 \pmod {17}$.
This~is enough to ensure that
$x_0/x_3 \in \bbQ_{17}^*$
is a~square. Since~the first equation shows that
$-x_0x_1$
is a norm from
$\smash{\bbQ_{17}({\!\sqrt{x_3/x_4}})}$,
one has
$\smash{(\frac{x_1}{x_4}, \frac{x_3}{x_4})_{17} = (-\frac{x_0}{x_4}, \frac{x_3}{x_4})_{17} = (-\frac{x_3}{x_4}, \frac{x_3}{x_4})_{17} = 1}$,
which completes the~proof.
\eop
\end{ex}

\begin{rems}
\label{intptex}
\begin{iii}
\item
Thus,~the Brauer class
$\tau$
causes a violation of strong approximation 
on~$X$
off~$\{17,\infty\}$.
In~fact, a
$\bbZ[\frac1{17}]$-valued
point such that
$x_1 \neq 0$
and
$x_3 \neq 0$
must necessarily fulfil
$(\frac{x_1}{x_4}, \frac{x_3}{x_4})_2 = 1$,
although not all
\mbox{$\bbZ_2$-valued}
points satisfy this~relation.
\item
This~restriction applies, in particular, to integral~points.
\item
There~are infinitely many integral points
on~$\calU$.
Indeed,~the curve defined by
$V(X_3)$
yields the~family
$$(-n^2\!:\!n^4\!:\!\pm n^3\!:\!0\!:\!1) \, .$$
On~the other hand, the curve defined by
$V(2X_1+X_2+X_3)$
is~elliptic. It~carries the six integral points
$(0\!:\!0\!:\!0\!:\!0\!:\!1)$,
$(0\!:\!0\!:\!-1\!:\!1\!:\!1)$,
$(-2\!:\!4\!:\!-1\!:\!-7\!:\!1)$,
$(-2\!:\!4\!:\!0\!:\!-8\!:\!1)$,
$(-14\!:\!196\!:\!-49\!:\!-343\!:\!1)$,
$(-14\!:\!196\!:\!48\!:\!-440\!:\!1)$,
and no~others.
\item
A rather naively implemented search for integral points 
on~$\calU$
delivered
$28$
of height
$<\! 50\,000$
that are not of the forms mentioned~above. These~are the~following,\newline
{\scriptsize
$(0\!:\!-2\!:\!1\!:\!1\!:\!1)$,
$(-1\!:\!-1\!:\!-1\!:\!2\!:\!1)$,
$(-1\!:\!-3\!:\!1\!:\!4\!:\!1)$,
$(4\!:\!-8\!:\!6\!:\!4\!:\!1)$,
$(-6\!:\!4\!:\!4\!:\!-8\!:\!1)$,
$(-8\!:\!18\!:\!-11\!:\!-23\!:\!1)$,
$(14\!:\!-28\!:\!20\!:\!8\!:\!1)$,
$(16\!:\!-56\!:\!30\!:\!4\!:\!1)$,
$(76\!:\!-696\!:\!-230\!:\!4\!:\!1)$,
$(-97\!:\!521\!:\!-223\!:\!-808\!:\!1)$,
$(-105\!:\!1413\!:\!381\!:\!-3204\!:\!1)$,
$(263\!:\!-829\!:\!467\!:\!62\!:\!1)$,
$(-556\!:\!912\!:\!712\!:\!-128\!:\!1)$,
$(-708\!:\!1278\!:\!-951\!:\!-423\!:\!1)$,
$(839\!:\!-1595\!:\!1157\!:\!444\!:\!1)$,
$(1004\!:\!-1648\!:\!-1288\!:\!4352\!:\!1)$,
$(-2073\!:\!3573\!:\!-2721\!:\!-2988\!:\!1)$,
$(2238\!:\!-6876\!:\!3924\!:\!9288\!:\!1)$,
$(-2264\!:\!3840\!:\!-2948\!:\!-3056\!:\!1)$,
$(-2916\!:\!5832\!:\!4122\!:\!-15228\!:\!1)$,
$(3324\!:\!-15678\!:\!-7219\!:\!289\!:\!1)$,
$(3879\!:\!-6183\!:\!-4899\!:\!16344\!:\!1)$,
$(-5450\!:\!14688\!:\!8947\!:\!-791\!:\!1)$,
$(-5809\!:\!11231\!:\!-8077\!:\!-2950\!:\!1)$,
$(8908\!:\!-16476\!:\!12115\!:\!5017\!:\!1)$,
$(-10194\!:\!6948\!:\!8415\!:\!-15687\!:\!1)$,
$(22238\!:\!-38044\!:\!29087\!:\!31097\!:\!1)$,
and
$(-26396\!:\!44152\!:\!-34138\!:\!-33148\!:\!1)$.}
\item
The exceptional behaviour at the prime
$17$
is possible only because
$D$
has bad reduction at the
prime~$17$.
The~curve
$D_{17}$
has exactly one singular point,
at~$(1\!:\!1\!:\!13\!:\!1\!:\!0)$.
Using Gr\"obner bases, one may easily verify that,
on~$D_{17}$,
there is the~relation
$$\frac{X_1}{X_3} = 8 \Big( \frac{X_0X_2X_3 - 2X_2^3 + X_0X_2^2 + X_0^2X_2}{X_0^2(4X_0+X_2)} \Big)^2 \, .$$
\end{iii}
\end{rems}

\begin{ex}
Let~$\calX' \subset \Pb^4_\bbZ$
be given by the system of~equations
\begin{align*}
                  X_0(8X_1+3X_4)+X_2^2 &= X_4(8X_3+2X_4) \, , \\
(8X_3+2X_4)(16X_1+X_2+8X_3+8X_4)+X_0^2 &= X_4(8X_1+3X_4) \, .
\end{align*}
Put
$\calU' := \calX' \!\setminus\! \calH$
for the hyperplane
$\calH := V(X_4) \subset \Pb^4_\bbZ$
and denote the generic fibre of
$\calU'$
by~$U'$.\smallskip

\noindent
Then~$\calU'(\bbZ_p) \neq \emptyset$
for every prime
number~$p$
and
$U'(\bbQ) \neq \emptyset$,
but there are no
\mbox{$\smash{\bbZ[\frac1{17}]}$-valued}
points,
$\smash{\calU'(\bbZ[\frac1{17}]) = \emptyset}$.
I.e,~the Hasse principle for
\mbox{$\smash{\bbZ[\frac1{17}]}$-valued}
points is~violated. In~particular, a failure of the integral Hasse principle~occurs. The~violations can be explained by a transcendental Brauer~class.\medskip

\noindent
{\bf Proof.}
There is a morphism
$\pi\colon \calU' \to \calU$
to the scheme from Example~\ref{ex1b} that is given by
$$\pi\colon (x_0\!:\!x_1\!:\!x_2\!:\!x_3\!:\!1) \mapsto (x_0\!:\!(8x_1+3)\!:\!x_2\!:\!(8x_3+2)\!:\!1) \, .$$
This~morphism yields bijections between the sets of
\mbox{$\bbQ$-rational}
points, as well as between the sets of
\mbox{$\bbZ_p$-valued}
points
for~$p \neq 2$.
In~particular,
$U'(\bbQ) \neq \emptyset$
and
$\calU'(\bbZ_p) \neq \emptyset$,
except possibly
for~$p \neq 2$.
Moreover,~there is the
\mbox{$\bbZ_2$-valued}
point
$(3\!:\!-\frac43\!:\!5\!:\!0\!:\!1) \in \calU'(\bbZ_2)$.

There~is, however, no
\mbox{$\smash{\bbZ[\frac1{17}]}$-valued}
point
on~$\calU'$.
In~fact, the image of such a point would fulfil
$(8x_1+3, 8x_3+2)_2 = (3, 2)_2 = -1$,
in contradiction to what was proven
about~$\calU$\!,
cf.\ Remark~\ref{intptex}.i).
\eop
\end{ex}

\begin{rems}
\label{notatinf}
\begin{iii}
\item
Observe that the non-existence of integral points
on~$\calU'$
may not be explained by a lack of real~points. In~fact, the
$\bbQ$-scheme
$U'$
is strongly unobstructed
at~$\infty$,
simply because
$U$~is.

\item
The result may be formulated in a more classical manner as~follows. The~pair
$(q_1, q_2)$
of quadratic polynomials in four variables, consisting of
\begin{align*}
q_1 &:= \hspace{3.67cm} X_0(8X_1+3) +X_2^2 - (8X_3+2) \, , \\
q_2 &:= (8X_3+2)(16X_1+X_2+8X_3+8) +X_0^2 - (8X_1+3) \, ,
\end{align*}
does not represent
$(0,0)$
integrally, although it does so over the rationals and in
\mbox{$p$-adic}
integers for every prime
number~$p$.
\end{iii}
\end{rems}

\subsubsection*{Blowing up a point. A cubic surface.}

\begin{ex}
\label{ex_cs1}
Let~$\calS \subset \Pb^3_\bbZ$
be given by the~equation
$$Y_0^3 + Y_0Y_1Y_2 + Y_0Y_2^2 - 2Y_1^2Y_2 + Y_1^2Y_3 + 2Y_2^2Y_3 - Y_2Y_3^2 = 0$$
and put
$\calV := \calS \!\setminus\! \calE$,
for the hyperplane
$\calE := V(Y_3) \subset \Pb^3_\bbZ$.\smallskip

\noindent
Then~every integral point
$(y_0\!:\!y_1\!:\!y_2\!:\!1) \in \calV(\bbZ)$
such that
$y_0y_2 \neq 0$
satisfies the condition that
$(y_0(y_2-y_1^2), y_2)_2 = 1$
or~$\gcd(y_0, 2y_2-1) > 1$.\medskip

\noindent
{\bf Proof.}
The scheme
$\calS$
is obtained from the scheme
$\calX$
of Example~\ref{ex1} by blowing up the point
$(0\!:\!1\!:\!0\!:\!0\!:\!0)$.
From~the computational viewpoint, this means to eliminate
$X_1$
from the equations
defining~$\calX$.
Moreover,~we replaced the coordinate sections
$X_0$,
$X_2$,
$X_3$,
$X_4$,
in this order, by
$Y_0$,
$Y_1$,
$Y_2$,
and~$Y_3$.
An~integral point
on~$\calV$
thus corresponds to a
\mbox{$\bbQ$-rational}
point
$(x_0\!:\!x_1\!:\!x_2\!:\!x_3\!:\!1) \in U(\bbQ)$
such that
$x_0$,
$x_2$,
and
$x_3$
are integers, but
$x_1$
not necessarily.

If,~however,
$\gcd(x_0, 2x_3-1) = 1$
then the equations
$x_0\!\cdot\!x_1 = x_3-x_2^2$
and
$(2x_3-1) \!\cdot\!x_1 = -(x_0^2 + x_2x_3 + x_3^2)$
together imply that
$x_1$
has to be an integer, as~well. Then,~by Example~\ref{ex1b}.i),
$(x_1, x_3)_2 = 1$.
The assertion follows, as
$\calS$
does not have integral points such that
$y_2-y_1^2 = 0$,
but~$y_0y_2 \neq 0$.
\eop
\end{ex}

\begin{rems}
\begin{iii}
\item
Write~$\smash{S \subset \Pb^3_\bbQ}$
and
$\smash{E \subset \Pb^3_\bbQ}$
for the generic fibres
of
$\calS$
and~$\calE$,
respectively, and put
$V := S \!\setminus\! E$.
Then,~under the blow-down morphism
$\pi\colon S \to X$,
one has
$\pi^{-1}(H) = E \cup L$
and, hence,
$\pi^{-1}(U) = V \!\setminus\! L$,
for
$L \subset S$
the line given by
$Y_0 = 2Y_2-Y_3 = 0$.
Note~that a point on the hyperplane
$H$
has been blown~up.

The effect deduced above by an elementary argument is, of course, that of the Brauer class
$\pi^* \tau \in \Br(V \!\setminus\! L)$.
The~points such that
$\gcd(y_0, 2y_2-1) = 1$
are exactly those that are integral with respect to the obvious integral model
of~$V \!\setminus\! L$.
\item
Blowing up a point outside
$H$
leads to
$\pi^{-1}(H)$
being a non-plane genus-one curve
on~$S$.
In~this case, the corresponding notion of integrality for points is further away from the classical~meaning.
\item
There~exist integral points
on~$\calV$
of all three kinds that are allowed by the statements~above. For~instance, for
$(-1\!:\!-1\!:\!2\!:\!1)$,
the
$\gcd$
is~$1$
and the Hilbert symbol
is~$1$.
On~the other hand,
for~$(-17\!:\!15\!:\!-8\!:\!1)$,
the
$\gcd$
is~$17>1$
and the Hilbert symbol
is~$1$,
while
for~$(3\!:\!5\!:\!2\!:\!1)$,
the
$\gcd$
is~$3>1$
and the Hilbert symbol is equal
to~$(-1)$.
\end{iii}
\end{rems}

\begin{ex}
Let~$\calS' \subset \Pb^3_\bbZ$
be given by the~equation
\begin{eqnarray*}
128Y_0^3 + 144Y_0^2Y_3 + 32Y_0Y_1Y_2 + 8Y_0Y_1Y_3 + 128Y_0Y_2^2 + 80Y_0Y_2Y_3 + 66Y_0Y_3^2 & \\
{} - 16Y_1^2Y_2- 3Y_1^2Y_3 - 4Y_1Y_2Y_3 + 80Y_2^2Y_3 + 40Y_2Y_3^2 + 12Y_3^3 &= 0
\end{eqnarray*}
and put
$\calV' := \calS' \!\setminus\! \calE$,
for the hyperplane
$\calE := V(Y_3) \subset \Pb^3_\bbZ$.\smallskip

\noindent
Then~every integral point
$(y_0\!:\!y_1\!:\!y_2\!:\!1) \in \calV'(\bbZ)$
satisfies
$\gcd(8y_0+3, 16y_2+3) > 1$.\medskip

\noindent
{\bf Proof.}
There is a morphism
$\pi\colon \calV' \to \calV$
to the scheme from Example~\ref{ex_cs1} that is given by
$$\pi\colon (y_0\!:\!y_1\!:\!y_2\!:\!1) \mapsto ((8y_0+3)\!:\!(2y_1+1)\!:\!(8y_2+2)\!:\!1) \, .$$
As~$(8y_0+3)[(8y_2+2)-(2y_1+1)^2] \equiv 3 \pmod 8$
and
$8y_2+2 \equiv 2 \pmod 8$,
the Hilbert symbol is always equal
to~$(-1)$,
as~required.
\eop
\end{ex}

\begin{rem}
The surface
$\calV'$
contains infinitely many integral points. Indeed, define the two sequences
$c$
and~$c'$
in~$\bbZ^3$
recursively~by
\begin{align}
\label{Pell}
\hspace{-0.8mm}
c_1 := [0,-2,0], \quad
c_2 := [-48,170,-24], \quad\hspace{3.3mm}
c_{i+2} := -110c_{i+1}-c_i-[48,48,24] \, , \nonumber \\
\hspace{-0.8mm}
c'_1 := [0,2,0], \quad\hspace{3.3mm}
c'_2 := [-48,-266,-24], \quad
c'_{i+2} := -110c'_{i+1}-c'_i-[48,48,24] \, .
\end{align}
Then, for each
$i \in \bbN$,
$(c_{i1}\!:\!c_{i2}\!:\!c_{i3}\!:\!1) \in \calV'(\bbZ)$
and
$(c'_{i1}\!:\!c'_{i2}\!:\!c'_{i3}\!:\!1) \in \calV'(\bbZ)$.
Observe~that the intersection of
$\calS'$
with the plane given by
$Y_0 = 2Y_2$
contains the exceptional line and, therefore, splits off a conic.

There are further integral points
on~$\calV'$
not being of this particular form, for instance
$(-1536\!:\!5414\!:\!-803\!:\!1)$
and~$(20\,706\!:\!-344\,632\!:\!534\!:\!1)$.
Moreover, both are the smallest members of infinite sequences of integral points of the same kind as~(\ref{Pell}). However,~the second member of the sequence starting at
$(-1536\!:\!5414\!:\!-803\!:\!1)$
involves 1340-digit integers, already.
\end{rem}

\section{More examples}
\label{sec_ex2}

\subsubsection*{A del Pezzo surface of degree four.}\leavevmode\\[1mm]
In the example below, an algebraic Brauer class interacts with transcendental ones.

\begin{ex}
\label{ex2}
Let~$X \subset \Pb^4_\bbQ$
be given by the system of~equations
\begin{align*}
      X_0X_1+X_2^2 &= X_4X_3 \, , \\
X_3(X_1+X_3)+X_0^2 &= X_4X_1
\end{align*}
and put
$U := X \!\setminus\! H$,
for the hyperplane
$H := V(X_4) \subset \Pb^4_\bbQ$.

\begin{iii}
\item
Then~$X$
is a del Pezzo surface of degree~four.
\item
$D := H \cap X$
is a non-singular curve of genus~one. Furthermore,
$$D(\bbQ) = \{(0\!:\!1\!:\!0\!:\!0\!:\!0), (0\!:\!1\!:\!0\!:\!-1\!:\!0)\} \, .$$
\item
The manifold\/
$X(\bbR)$
is connected, its submanifold
$U(\bbR)$
consists of two connected components, and
$U$
is strongly unobstructed
at~$\infty$.
\item
One has
$\Br(U)/\Br_0(U) \cong (\bbZ/2\bbZ)^2$,
the three nontrivial elements being represented~by
\begin{align*}
       \alpha\colon &\textstyle \;(\frac{X_1}{X_4}, -1; -1) \, , \\
         \tau\colon &\textstyle \;(\frac{X_1}{X_4}, \frac{X_3}{X_4}; -1) \, , \qquad {\rm and} \\
  \alpha+\tau\colon &\textstyle \;(\frac{X_1}{X_4}, -\frac{X_3}{X_4}; -1) \, .
\end{align*}
Moreover,~$\Br_1(U)/\Br_0(U) = \langle\alpha\rangle \cong \bbZ/2\bbZ$.
\end{iii}\medskip

\noindent
{\bf Proof.}
i) and~ii)
{\tt magma} reports that both
$X$
and~$D$
are non-singular. More~precisely,
$D$~is
isomorphic to the elliptic curve
$E\colon y^2z = x^3 + xz^2$
of
$j$-invariant
$1728$,
an isomorphism
$\iota\colon D \to E$
being given by
$$(X_0\!:\!X_1\!:\!X_2\!:\!X_3\!:\!0) \mapsto (X_2^2 \!:\! -X_1X_2 \!:\! -X_1X_3) \, .$$
This~explains, in particular, why
$D$
has no
\mbox{$\bbQ$-rational}
points other than the two obvious~ones.\smallskip

\noindent
iii)
The~pencil of quadrics
in~$\Pb^4$
corresponding
to~$X$
contains exactly three degenerate quadrics that are defined
over~$\bbR$.
Among~the corresponding rank-four discriminants, two are negative but the third one is~positive. According to~\cite[Theorem~3.4 and Lemma~5.1]{VAV}, this already implies that
$\Br(X_\bbR)/\Br(\bbR) = 0$.
Consequently,~$X(\bbR)$
is connected~\cite[Th\'eor\`eme 1.4]{Si}.

As~the cubic polynomial
$x^3+x$
has only one real root, we see that
$D(\bbR)$
is connected,~too. Such~a curve may cut
$X(\bbR)$
into not more than two~components. On~the other
hand,~$U(\bbR)$
is clearly~disconnected. Indeed,~for a real point
$(x_0\!:\!x_1\!:\!x_2\!:\!x_3\!:\!1) \in U(\bbR)$,
the second equation implies
$(x_3 + \frac{x_1}2)^2 + x_0^2 = \frac14[x_1^2 + 4x_1]$,
hence
$x_1 \geq 0$
or~$x_1 \leq -4$.
And~there exist points of both kinds, for example
$(-1\!:\!1\!:\!0\!:\!-1\!:\!1)$
and~$(0\!:\!-4\!:\!\sqrt{2}\!:\!2\!:\!1)$.

Since~$D$
is a non-singular, geometrically irreducible curve having real points and both components of
$U(\bbR)$
have limit points
in~$D(\bbR)$,
Theorem~\ref{unobstr} yields that
$U$~is
strongly unobstructed at the real~place.\smallskip

\noindent
iv)
A~calculation using Gr\"obner bases shows that the Galois group operating on the 16~lines
on~$X$
is of order
exactly~$384$.
Thus,~Theorem \ref{algBr}.b.ii) immediately implies that
$\Br_1(U)/\Br_0(U) \cong \bbZ/2\bbZ$.
Moreover,~the second equation is only of
rank~$4$
and may be written in the form
$X_1(X_4-X_3) - X_3^2 - X_0^2 = 0$.
Therefore,~Theorem~\ref{algBrexpl} shows that
$\smash{(\frac{X_1}{X_4}, -1; -1)}$
defines a nontrivial algebraic Brauer class
on~$U$.

On~the other hand, we have
$\#D(\bbQ) = 2$.
Again,~$E$
has two nontrivial
\mbox{$4$-isogenies,}
but no nontrivial
\mbox{$n$-isogenies}
for any
other~$n \geq 3$.
Only~one of the isogenous curves has a
\mbox{$\bbQ$-rational}
\mbox{$4$-torsion}
point, which, however, is not in the image of the
\mbox{$4$-torsion}
of~$E$.
Hence,~Corollary~\ref{trBrcor} provides us with an injection
$\Br(U)/\Br_1(U) \hookrightarrow \bbZ/2\bbZ$.
Consequently,~$\Br(U)/\Br_0(U)$
is of order at
most~$4$.
As~far as a concrete description is asked for, Theorem~\ref{trBrexpl} applies. It~shows that
$\smash{(\frac{X_1}{X_4}, \frac{X_3}{X_4}; -1)}$
defines a Brauer class
on~$U$.

In~order to prove that this is indeed a transcendental class, we need to show that
$\smash{\frac{X_1}{X_3}}$
is not the square of a rational function
on~$D_{\overline\bbQ}$.
For~this, we may argue as~follows.
One~has
\begin{align*}
\textstyle -\frac{X_1}{X_3} &=\textstyle -\frac{X_1X_3}{X_3^2} \\
                            &=\textstyle \frac{X_0^2+X_3^2}{X_3^2}  \\
                            &=\textstyle \frac{X_0^2X_1^2+X_1^2X_3^2}{X_1^2X_3^2} \\
                            &=\textstyle \frac{X_2^4+X_1^2X_3^2}{X_1^2X_3^2}
\end{align*}
and this is the product of
$\smash{F := \frac{X_2^2+iX_1X_3}{X_1X_3}}$
together with its complex~conjugate.
One~directly finds that
$$\div F = 2[(1\!:\!0\!:\!0\!:\!i\!:\!1)] - 2[(0\!:\!1\!:\!0\!:\!0\!:\!1)] \, ,$$
hence
$$\div \overline{F} = 2[(1\!:\!0\!:\!0\!:\!-i\!:\!1)] - 2[(0\!:\!1\!:\!0\!:\!0\!:\!1)] \, ,$$
and
$$\textstyle \div (\frac{X_1}{X_3}) = 2[(1\!:\!0\!:\!0\!:\!i\!:\!1)] + 2[(1\!:\!0\!:\!0\!:\!-i\!:\!1)] - 4[(0\!:\!1\!:\!0\!:\!0\!:\!1)] \, .$$
In~particular,
$(1\!:\!0\!:\!0\!:\!i\!:\!1) - (0\!:\!1\!:\!0\!:\!0\!:\!1) \in J(D)(\overline\bbQ)$
is a proper
\mbox{$2$-torsion}
point, and
$\smash{(1\!:\!0\!:\!0\!:\!-i\!:\!1) - (0\!:\!1\!:\!0\!:\!0\!:\!1) \in J(D)(\overline\bbQ)}$
is~another.

Therefore,~their sum
$(1\!:\!0\!:\!0\!:\!i\!:\!1) + (1\!:\!0\!:\!0\!:\!-i\!:\!1) - 2(0\!:\!1\!:\!0\!:\!0\!:\!1) \in J(D)(\overline\bbQ)$
is the third
\mbox{$2$-torsion}
point, in particular it is non-zero. As~a consequence of this, we see that
$\smash{\frac{X_1}{X_3}}$
can not be the square of a rational function
on~$D_{\overline\bbQ}$,
as~required.
\eop
\end{ex}

\begin{ex}[continued]
\label{ex2b}
Let~$X$,
$U$,
$\alpha$,
and~$\tau$
be as~above. Moreover,~let
$\calX \subset \Pb^4_\bbZ$
be defined by the same system of equations
as~$X$
and put
$\calU := \calX \!\setminus\! \calH$,
for
$\calH := V(X_4) \subset \Pb^4_\bbZ$.
Then

\begin{iii}
\item
for all primes
$p \neq 2, \infty$,
the local evaluation maps
$\ev_{\alpha, \nu_p}\colon \calU(\bbZ_p) \to \bbQ/\bbZ$,
$\ev_{\tau, \nu_p}\colon \calU(\bbZ_p) \to \bbQ/\bbZ$,
and
$\ev_{\alpha+\tau, \nu_p}\colon \calU(\bbZ_p) \to \bbQ/\bbZ$
are constantly~zero.
\item
The local evaluation map
$\ev_{\tau, \nu_\infty}\colon U(\bbR) \to \frac12\bbZ/\bbZ$
is constantly zero, but
$\ev_{\tau, \nu_2}\colon \calU(\bbZ_2) \to \bbQ/\bbZ$
is non-constant.
\item
The local evaluation map
$\ev_{\alpha+\tau, \nu_\infty}\colon U(\bbR) \to \frac12\bbZ/\bbZ$
is non-constant, but
$\ev_{\alpha+\tau, \nu_2}\colon \calU(\bbZ_2) \to \bbQ/\bbZ$
is constantly zero.
\item
The maps
$\ev_{\alpha, \nu_2}\colon \calU(\bbZ_2) \to \bbQ/\bbZ$
and
$\ev_{\alpha, \nu_\infty}\colon U(\bbR) \to \frac12\bbZ/\bbZ$
are both non-constant.
\item
For~$p \equiv 1 \pmod 4$,
even the local evaluation map
$\ev_{\alpha, \nu_p}\colon U(\bbQ_p) \to \bbQ/\bbZ$
is constantly~zero.
\end{iii}\medskip

\noindent
{\bf Proof.}
i)
For~$\ev_{\tau, \nu_p}$,
this is directly Lemma~\ref{trconst}. And for
$\ev_{\alpha, \nu_p}$,
Lemma~\ref{algconst} applies, since
$d = -1$,
the linear forms
$X_1$,
$X_4-X_3$,
$X_3$,
and~$X_0$
are linearly independent modulo every
prime~$p$,
and the cusp
$(0\!:\!0\!:\!1\!:\!0\!:\!0)$
does not satisfy the second equation
$X_0X_1 + X_2^2 = X_4X_3$
modulo any~prime.\smallskip

\noindent
ii)
The~second equation excludes the possibility that
$X_1$
and
$X_3$
might both be~negative.
Hence,~$\ev_{\tau, \nu_\infty}$
is constantly~zero.

Furthermore,~for the integral point
$x = (-1\!:\!1\!:\!0\!:\!-1\!:\!1) \in \calU(\bbZ_2)$,
one has
$(x_1, x_3)_2 = (1, -1)_2 = 1$
and thus
$\ev_{\tau, \nu_2}(x) = 0$.
On~the other hand, the
\mbox{$\bbZ_2$-valued}
point
$\smash{x' = (-\frac13\!:\!\frac13\!:\!-\frac23\!:\!\frac13\!:\!1) \in \calU(\bbZ_2)}$
yields
$(x_1, x_3)_2 = (\frac13, \frac13)_2 = -1$,
hence
$\ev_{\tau, \nu_2}(x) = \frac12$.\smallskip

\noindent
iii)
First~of all, the integral point
$x = (-1\!:\!1\!:\!0\!:\!-1\!:\!1) \in U(\bbR)$
yields
$(x_1, -x_3)_\infty = (1, 1)_\infty = 1$,
therefore
$\ev_{\alpha+\tau, \nu_\infty}(x) = 0$.
On~the other hand, there is the real point
$\smash{x' = (0\!:\!-4\!:\!\sqrt{2}\!:\!2\!:\!1)}$
such that
$(x'_1, -x'_3)_\infty = (-4, -2)_\infty = -1$,
hence
$\smash{\ev_{\alpha+\tau, \nu_\infty}(x) = \frac12}$.

Finally,~constancy of
$\ev_{\alpha+\tau, \nu_2}$
is elementary, too, but quite~involved. The~idea of the proof is roughly as~follows. It~suffices to show
$(x_1, -x_3)_2 = 1$
for every
\mbox{$\bbZ_2$-valued}
point
$\smash{x = (x_0\!:\!x_1\!:\!x_2\!:\!x_3\!:\!x_4) \in \calU(\bbZ_2)}$
such that
$x_1 \neq 0$
and~$x_3 \neq 0$.

In~order to do this, several cases have to be distinguished. It~turns out that the assumption
$8 \notd x_3$
leads to only finite many cases and that
$(x_1, -x_3)_2 = 1$
is fulfilled in each of~them. On~the other hand,
if~$8 | x_3$
then the second equation
$x_0^2 + x_3^2 = (1-x_3)x_1$
shows that
$(x_1, -1)_2 = 1$.
Moreover,~$x_1$
is automatically a square, except for the case that
$|\nu_2(x_0) - \nu_2(x_3)| \leq 1$.
In~this case, however,
$\nu_2(x_0x_1) \geq \nu_2(x_3)+3$
such that the first equation implies that
$x_3$
is a square. The~assertion~follows.\smallskip

\noindent
iv)
is an immediate consequence of ii) and~iii), while v) follows from the fact that
$(-1)$
is a square
in~$\bbQ_p$
for~$p \equiv 1 \pmod 4$.
\eop
\end{ex}

\begin{rems}
\label{intptex2}
\begin{iii}
\item
The algebraic Brauer class
$\alpha$
causes a violation of strong approximation off
$S_1 := \{ p \text{~prime} \mid p \equiv 1 \pmod 4 \}$.
In fact, a
\mbox{$\smash{\bbZ[\frac1{S_1}]}$-valued}
point such that
$x_1 \neq 0$
must necessarily fulfil
$(\frac{x_1}{x_4},-1)_2 + (\frac{x_1}{x_4},-1)_\infty = 0$,
although not all adelic points outside
$S_1$
satisfy this~relation.
\item
Furthermore,~strong approximation off
$\{\infty\}$
is violated. In~fact, the Brauer
class~$\tau$
yields a restriction for integral~points. Namely,~an integral point such that
$x_1 \neq 0$
and
$x_3 \neq 0$
must necessarily fulfil
$(x_1, x_3)_2 = 1$,
although not all
\mbox{$\bbZ_2$-valued}
points satisfy this~relation.
\item
The Brauer class
$\alpha+\tau$
yields another restriction for integral~points. Indeed,~an integral point such that
$x_1 \neq 0$
and
$x_3 \neq 0$
must necessarily fulfil
$(x_1, -x_3)_\infty = 1$,
although there exist real points violating~this.
\item
There~are infinitely many integral points
on~$\calU$.
In~fact, the curves defined by
$V(X_3)$
and
$V(X_1+X_3)$
yield the~families
$$(-n^2\!:\!n^4\!:\!\pm n^3\!:\!0\!:\!1) \quad{\rm and}\quad (-(n^2+1)\!:\!(n^2+1)^2\!:\!\pm (n^2+1)n\!:\!-(n^2+1)^2\!:\!1) \, .$$
\item
A~search for integral points 
on~$\calU$
delivered exactly 16 of height
$<\! 50\,000$
that are of neither of the two types above. These~are\newline
{\scriptsize
$(-5\!:\!13\!:\!\pm8\!:\!-1\!:\!1)$,
$(-58\!:\!676\!:\!\pm198\!:\!-4\!:\!1)$,
$(-268\!:\!4240\!:\!\pm1064\!:\!-4224\!:\!1)$,
$(-1297\!:\!11437\!:\!\pm3850\!:\!-11289\!:\!1)$,
$(-2416\!:\!6736\!:\!\pm4034\!:\!-1020\!:\!1)$,
$(-4513\!:\!9685\!:\!\pm6611\!:\!-3084\!:\!1)$,
$(-6668\!:\!13456\!:\!\pm9472\!:\!-5824\!:\!1)$,
and
$(-11681\!:\!27061\!:\!\pm17779\!:\!-6700\!:\!1)$.}
\end{iii}
\end{rems}

\begin{ex}
Let~$\calX' \subset \Pb^4_\bbZ$
be given by the system of~equations
\begin{align*}
     X_0(4X_1+3X_4)+X_2^2 &= X_4(4X_3+3X_4) \, , \\
X_3(4X_1+4X_3+6X_4)+X_0^2 &= X_4(4X_1+3X_4)
\end{align*}
and put
$\calU' := \calX' \!\setminus\! \calH$,
for the hyperplane
$\calH := V(X_4) \subset \Pb^4_\bbZ$.
Denote~the generic fibre of
$\calU'$
by~$U'$.\smallskip

\noindent
Then~$\calU'(\bbZ_p) \neq \emptyset$
for every prime
number~$p$
and
$U'(\bbQ) \neq \emptyset$,
but there are no integral~points,
$\calU'(\bbZ) = \emptyset$.
In~other words, the integral Hasse principle is~violated. The~violation can be explained by a Brauer-Manin~obstruction.\medskip

\noindent
{\bf Proof.}
There is a morphism
$\pi\colon \calU' \to \calU$
to the scheme from the example before that is given by
$$\pi\colon (x_0\!:\!x_1\!:\!x_2\!:\!x_3\!:\!1) \mapsto (x_0\!:\!(4x_1+3)\!:\!x_2\!:\!(4x_3+3)\!:\!1) \, .$$
This~morphism yields bijections between the sets of
\mbox{$\bbQ$-rational}
points, as well as between the sets of
\mbox{$\bbZ_p$-valued}
points
for~$p \neq 2$.
In~particular,
$U'(\bbQ) \neq \emptyset$
and
$\calU'(\bbZ_p) \neq \emptyset$,
except possibly
for~$p \neq 2$.
Moreover,~there is a
\mbox{$\bbZ_2$-valued}
point such that
$x_1 = 0$
and~$x_3 = 2$.

There~is, however, no integral point
on~$\calU'$.
In~fact, the image of such a point would fulfil
$(4x_1+3, 4x_3+3)_2 = (3,3)_2 = -1$,
in contradiction to what was proven, cf.\ Remark~\ref{intptex2}.i).
\eop
\end{ex}

\begin{rem}
Again,~the
\mbox{$\bbQ$-scheme}
$U'$
is strongly unobstructed
at~$\infty$.
Concerning~this, the situation is completely analogous to the one described in Remark~\ref{notatinf}.i).
\end{rem}

\subsubsection*{Blowing up a point. A cubic surface.}

\begin{ex}
\label{ex_cs2}
Let~$\calS \subset \Pb^3_\bbZ$
be given by the~equation
$$Y_0^3 + Y_0Y_2^2 - Y_1^2Y_2 + Y_1^2Y_3 + Y_2^2Y_3 - Y_2Y_3^2 = 0 \, .$$
Put
$\calV := \calS \!\setminus\! \calE$
for the hyperplane
$\calE := V(Y_3) \subset \Pb^3_\bbZ$.\smallskip

\noindent
Then~every integral point
$(y_0\!:\!y_1\!:\!y_2\!:\!1) \in \calV(\bbZ)$
such that
$y_0y_2 \neq 0$
satisfies that
$(y_0(y_2-y_1^2), y_2)_2 = 1$.
In~particular, strong approximation
off~$\{\infty\}$
is~violated.\medskip

\noindent
{\bf Proof.}
The scheme
$\calS$
is obtained from the scheme
$\calX$
from Example~\ref{ex2} by blowing up the point
$(0\!:\!1\!:\!0\!:\!0\!:\!0)$.
We~eliminated
$X_1$
from the equations
defining~$\calX$
and replaced the coordinate sections
$X_0$,
$X_2$,
$X_3$,
$X_4$,
in this order, by
$Y_0$,
$Y_1$,
$Y_2$,
and~$Y_3$.
An~integral point
on~$\calV$
corresponds to a
\mbox{$\bbQ$-rational}
point
$(x_0\!:\!x_1\!:\!x_2\!:\!x_3\!:\!1) \in U(\bbQ)$
such that
$x_0$,
$x_2$,
and
$x_3$
are integers, but
$x_1$
not necessarily.

For~the main assertion, let us first observe that there are no integral points such that
$y_2-y_1^2 = 0$,
but~$y_0y_2 \neq 0$.
Thus,~the symbol is properly defined. Moreover, Example~\ref{ex2b}.ii) shows that
$(y_0(y_2-y_1^2), y_2)_\infty = 1$.
We~will complete the proof by verifying
$(y_0(y_2-y_1^2), y_2)_p = 1$
for every
prime~$p \neq 2$.

In~order to do this, assume that
$p \notd (x_3-1)$~first.
Then~$(x_3-1) \!\cdot\!x_1 = -(x_0^2 + x_3^2)$
implies that
$x_1$
is a
\mbox{$p$-adic}
integer. Consequently, by Example~\ref{ex2b}.i), one has
$(y_0(y_2-y_1^2), y_2)_p = (x_1, x_3)_p = 1$.
On~the other hand, if
$p \,|\, (x_3-1)$
then
$x_3$
is a square
in~$\bbZ_p^*$
and therefore
$(y_0(y_2-y_1^2), y_2)_p = (x_1, x_3)_p = 1$,
as~well. The assertion~follows.

Finally,~we note that, for the point
$\smash{(y_0\!:\!y_1\!:\!y_2\!:\!1) = (-\frac13\!:\!-\frac23\!:\!\frac13\!:\!1) \in \calV(\bbZ_2)}$,
one has that
$(y_0(y_2-y_1^2), y_2)_2 = -1$.
Hence,~this
\mbox{$\bbZ_2$-valued}
point cannot be approximated by integral~ones.
\eop
\end{ex}

\begin{rems}
Write~$\smash{S \subset \Pb^3_\bbQ}$
and
$\smash{E \subset \Pb^3_\bbQ}$
for the generic fibres
of
$\calS$
and~$\calE$,
respectively, and put
$V := S \!\setminus\! E$.
Then,~under the blow-down morphism
$\pi\colon S \to X$,
one has
$\pi^{-1}(H) = E \cup L$
and, hence,
$\pi^{-1}(U) = V \!\setminus\! L$,
for
$L \subset S$
the line given by
$Y_0 = Y_2-Y_3 = 0$.

\begin{iii}
\item
The~effect shown above is that of the Brauer class
$\pi^* \tau \in \Br(V \!\setminus\! L)$,
which, as
$\smash{\frac{Y_2}{Y_3} = \frac{X_3}{X_4}}$
is a square
on~$L$,
turns out to be unobstructed
at~$L$
and extends to the whole
of~$V$.
\item
The other generator shows a behaviour similar to that indicated in Example~\ref{ex_cs2}. I.e.,~it yields that
$\gcd(y_0, y_2-1) > 1$
or~$(y_0(y_2-y_1^2), -y_2)_\infty = 1$
for every integral point, while neither of the two statements is always~true.
\end{iii}
\end{rems}

\begin{ex}
Let~$\calS' \subset \Pb^3_\bbZ$
be given by the~equation
\begin{eqnarray*}
16Y_0^3 + 12Y_0^2Y_3 + 16Y_0Y_2^2 + 24Y_0Y_2Y_3 + 12Y_0Y_3^2 - 4Y_1^2Y_2 - 2Y_1^2Y_3 + 8Y_2^2Y_3 & \\
{} + 11Y_2Y_3^2 + 4Y_3^3 &= 0
\end{eqnarray*}
and put
$\calV' := \calS' \!\setminus\! \calE$,
for the hyperplane
$\calE := V(Y_3) \subset \Pb^3_\bbZ$.\smallskip

\noindent
Then~$\calV'(\bbZ_p) \neq \emptyset$
for every prime
number~$p$
and
$\calV'(\bbQ) \neq \emptyset$,
but there are no integral~points,
$\calV'(\bbZ) = \emptyset$.
In~other words, the integral Hasse principle is~violated. The~violation is explained by a Brauer-Manin~obstruction.\medskip

\noindent
{\bf Proof.}
There is a morphism
$\pi\colon \calV' \to \calV$
to the scheme from Example~\ref{ex_cs2} that is given by
$$\pi\colon (y_0\!:\!y_1\!:\!y_2\!:\!1) \mapsto ((4y_0+1)\!:\!2y_1\!:\!(4y_2+3)\!:\!1) \, .$$
This~morphism yields bijections between the sets of
\mbox{$\bbQ$-rational}
points, as well as between the sets of
\mbox{$\bbZ_p$-valued}
points
for~$p \neq 2$.
In~particular,
$\calV'(\bbQ) \neq \emptyset$
and
$\calV'(\bbZ_p) \neq \emptyset$,
except possibly
for~$p \neq 2$.
Moreover,~there is the
\mbox{$\bbZ_2$-valued}
point
$(-\frac13\!:\!-\frac13\!:\!-\frac23\!:\!1)$.

There~is, however, no integral point
on~$\calV'$.
In~fact, the image of such a point would fulfil
$((4y_0+1)[4y_2+3 - 4y_1^2], 4y_2+3)_2 = (3,3)_2 = -1$,
in contradiction to what was~proven.
\eop
\end{ex}

\section{Yet another example}
\label{sec_ex3}

In the example below, algebraic Brauer classes interact with effects caused by
$U(\bbR)$
being disconnected into compact and non-compact~components.
The two non-compact components
of~$U(\bbR)$
in fact fulfil the requirements of Definition~\ref{def_obst}.ii). Nevertheless,~strong approximation
off~$\{\infty\}$
is violated, as there are Brauer classes
$\alpha_1$
and~$\alpha_2$
working at the primes
$2$,
$3$,
and~$\infty$.

Just~requiring
\mbox{$2$-adic}
and
\mbox{$3$-adic}
approximation by a sequence
$(x_n)_{n \in \bbN}$
of integral points, these classes enforce a certain behaviour at the infinite prime, as~well. In~the present case, they determine the connected component
of~$U(\bbR)$
the points
$x_n$
are lying~on. The~violation of strong approximation
off~$\{\infty\}$
happens when trying to approximate an adelic point, for which the
$x_n \in \calU(\bbZ)$,
for~$n\to\infty$,
turn out to be bound to to the single compact component
of~$U(\bbR)$.
This~underlines observations made by U.~Derenthal and D.~Wei in \cite[Example~6.2]{DW}.

As~a consequence, we conclude that every serious definition of unobstructedness at infinity must include requirements on all connected components
of~$U(\bbR)$.
In~particular, it is clearly insufficient just to require the existence of a single well-behaved~one. Or~to require that each irreducible component of the boundary contains limit points just
from~$U(\bbR)$.

\begin{ex}
Let~$X \subset \Pb^4_\bbQ$
be given by the system of~equations
\begin{align*}
        X_0(X_0+X_1) &= X_2^2 + (X_0+X_4)^2 \, , \\
 (X_0+X_2)(X_0+2X_2) &= 2X_1^2 + 3X_3^2
\end{align*}
and put
$U := X \!\setminus\! H$,
for the hyperplane
$H := V(X_4) \subset \Pb^4_\bbQ$.

\begin{iii}
\item
Then~$X$
is a del Pezzo surface of degree~four.
\item
The~curve
$D := H \cap X$
is non-singular of genus~one. One~has
$D(\bbQ) = \emptyset$,
but
$D(\bbR) \neq \emptyset$.
Furthermore,
$J(D)(\bbQ) \cong \bbZ$.
\item
The manifold\/
$X(\bbR)$
consists of two connected components and its submanifold
$U(\bbR)$
decomposes into three connected components.
$U$~is
not unobstructed
at~$\infty$,
not even in the weak~sense.
\item
One has
$\Br(U)/\Br_0(U) \cong (\bbZ/2\bbZ)^2$,
two generators being represented~by
\begin{align*}
  \alpha_1\colon &\textstyle \;(\frac{X_0}{X_4}, -1; -1) \qquad {\rm and} \\
  \alpha_2\colon &\textstyle \;(\frac{X_0+X_2}{X_4}, -6; -1) \, .
\end{align*}
Moreover,~$\Br(U) = \Br_1(U)$.
\end{iii}\medskip

\noindent
{\bf Proof.}
i) and~ii)
{\tt magma} reports that both
$X$
and~$D$
are non-singular. Furthermore,~the system of equations defining
$D$
is easily seen to be insoluble
in~$\bbQ_3$.
Therefore~$\smash{D(\bbQ) = \emptyset}$.
On~the other hand,
$\smash{(1\!:\!1\!:\!1\!:\!\frac23\sqrt{3}\!:\!0) \in D(\bbR)}$.
The~Jacobian
of~$D$
is reported to be isomorphic to the elliptic curve
$\smash{E\colon y^2z = x^3 + 6x^2z + 72xz^2 - 54z^3}$
over~$\bbQ$
of
$j$-invariant
$\frac{512\,000}{603}$
and Mordell-Weil group
$E(\bbQ)$
being free of rank~one.\smallskip

\noindent
iii)
The~pencil of quadrics
in~$\Pb^4$
defining~$X$
contains five degenerate ones, all of which are~real. Among~the corresponding rank-four discriminants, four are negative and the fifth is~positive. According to~\cite[Theorem~3.4 and Lemma~5.1]{VAV}, this implies
$\Br(X_\bbR)/\Br(\bbR) = (\bbZ/2\bbZ)^2$,
which in itself shows that
$X(\bbR)$
consists of exactly two connected components~\cite[Th\'eor\`eme 1.4]{Si}. The~evaluation of the global Brauer class
on~$X_\bbR$,
given by
$\smash{(\frac{X_0+X_2}{X_0}, -1; -1)}$,
distinguishes the two~components.

As~the cubic polynomial
$x^3 + 6x^2 + 72x - 54$
has exactly one real root, we see that
$D(\bbR)$
is connected. Such~a curve may meet only one component of
$X(\bbR)$,
it obviously meets
$X_+$,
the one where
$\smash{\frac{X_0+X_2}{X_0}}$
is positive, and it cannot cut that into more than two~components. On~the other hand, one sees that
$D(\bbR)$
indeed decomposes
$X_+$
into two components, these being distinguished by the sign of
$\smash{\frac{X_0}{X_4}}$.
Observe that both are non-empty, as they contain the real points
$(17\!:\!3\!:\!4\!:\!-13\!:\!1)$
and
$\smash{(-\frac{41}8\!:\!\frac32\!:\!\frac54\!:\!-\frac{11}8\!:\!1)}$,
respectively.

Finally,~let us note once again that the component
$X_- := X(\bbR) \!\setminus\! X_+$
of~$X(\bbR)$
is not at all met by the
curve~$D(\bbR)$.
Hence,~$X_{-}$
forms as well a connected component
of~$U(\bbR)$.
As~$X_-$
is compact, Definition~\ref{def_obst}.i) immediately shows that
$U$
cannot be weakly~unobstructed.\smallskip

\noindent
iv)
A~calculation using Gr\"obner bases indicates that the Galois group operating on the 16~lines
on~$X$
is of order
exactly~$96$.
As~both quadrics used to
describe~$X$
are defined
over~$\bbQ$
and of
rank~$4$,
Theorem \ref{algBr}.b.iii) implies that
$\Br_1(U)/\Br_0(U) \cong (\bbZ/2\bbZ)^2$.

Moreover,~Theorem~\ref{algBrexpl} shows that
$\smash{(\frac{X_0}{X_4}, -1; -1)}$
and
$\smash{(\frac{X_0+X_2}{X_4}, -6; -1)}$
define nontrivial algebraic Brauer classes
on~$U$.
They~do not just differ by an element
from~$\Br(\bbQ)$,
as easily follows from the results given in Example~\ref{ex3b}.i) and~ii), below. On~the other hand,
$J(D)$
do not have any nontrivial
\mbox{$\bbQ$-rational}
isogeny. Thus,~Corollary~\ref{trBrcor} yields
$\Br(U) = \Br_1(U)$.
\eop
\end{ex}

\begin{rems}
\begin{iii}
\item
Writing~down the points on
$U$
in the form
$(x_0\!:\!x_1\!:\!x_2\!:\!x_3\!:\!1)$,
the three connected components
of~$U(\bbR)$
are distinguished by the signs of
$x_0$
and
$x_0+x_2$
(respectively that of
$x_0+2x_2$
if~$x_0+x_2 = 0$).
\item
There are no real points
on~$U$
such that
$x_0 > 0$
and
$x_0+x_2 < 0$,
in accordance with our result that
$U(\bbR)$
consists of only three~components. One~may see this in an elementary fashion, as~follows.

The~first equation implies
$\smash{x_1 = \frac{x_2^2}{x_0} + \frac1{x_0} + 2}$,
hence
$\smash{x_1 > \frac{x_2^2}{x_0}}$
in the case
that~$x_0 > 0$.
Moreover,~$\smash{(\frac{x_0}{x_2})^2 (\frac{x_0}{x_2}+1) (\frac{x_0}{x_2}+2) > 2}$,
which implies that
$\smash{\frac{x_0}{x_2} < -2.293\ldots < -1}$
if
$x_2$
is negative, and hence
$x_0+x_2 = x_2(\frac{x_0}{x_2} + 1) > 0$.
\item
One may as well verify in an elementary manner that the component where
$x_0 < 0$
and
$x_0+x_2 > 0$
is~compact. All~such real points fulfil
$x_2 < 3$.

Indeed,~the equations
of~$X$
imply that
$\smash{(\frac{x_0}{x_2})^2 (\frac{x_0}{x_2}+1) (\frac{x_0}{x_2}+2) \geq 2[1+\frac1{x_2^2}+\frac{2x_0}{x_2^2}]^2}$.
If one assumes
$x_2 \geq 3$
and
$|x_0| < |x_2|$
then the right hand side is strictly larger than
$2[1+0-\frac23]^2 = \frac29$.
However,~the polynomial
$t^2(t+1)(t+2)$
adopts its maximum in the
interval~$[-1,0]$
at the point
$\smash{t = \frac{-9+\sqrt{17}}8}$
and this maximum is
$0.201\ldots < \frac29$.
\end{iii}
\end{rems}

\begin{ex}[continued]
\label{ex3b}
Let~$X$,
$U$,
$\alpha_1$,
and~$\alpha_2$
be as~above. Moreover,~let
$\calX \subset \Pb^4_\bbZ$
be defined by the same system of equations
as~$X$
and put
$\calU := \calX \!\setminus\! \calH$,
for
$\calH := V(X_4) \subset \Pb^4_\bbZ$.
Then

\begin{iii}
\item
for all primes
$p \neq 2, 3, \infty$,
the local evaluation maps
$\ev_{\alpha_1, \nu_p}\colon \calU(\bbZ_p) \to \bbQ/\bbZ$,
$\ev_{\alpha_2, \nu_p}\colon \calU(\bbZ_p) \to \bbQ/\bbZ$,
and
$\ev_{\alpha_1+\alpha_2, \nu_p}\colon \calU(\bbZ_p) \to \bbQ/\bbZ$
are constantly~zero.
\item
The local evaluation map
$\ev_{\alpha_1, \nu_3}\colon \calU(\bbZ_3) \to \bbQ/\bbZ$
is constantly zero, while
$\smash{\ev_{\alpha_1, \nu_\infty}\colon U(\bbR) \to \frac12\bbZ/\bbZ}$
and
$\ev_{\alpha_1, \nu_2}\colon \calU(\bbZ_2) \to \bbQ/\bbZ$
are non-constant.
\item
The local evaluation maps
$\smash{\ev_{\alpha_2, \nu_\infty}\colon U(\bbR) \to \frac12\bbZ/\bbZ}$,
$\ev_{\alpha_2, \nu_2}\colon \calU(\bbZ_2) \to \bbQ/\bbZ$,
and
$\ev_{\alpha_2, \nu_3}\colon \calU(\bbZ_3) \to \bbQ/\bbZ$
are all non-constant.
\item
The local evaluation maps
$\ev_{\alpha_1+\alpha_2, \nu_\infty}\!\colon U(\bbR) \to \!\frac12\bbZ/\bbZ$,
$\ev_{\alpha_1+\alpha_2, \nu_2}\!\colon \calU(\bbZ_2) \to \!\bbQ/\bbZ$,
and
$\ev_{\alpha_1+\alpha_2, \nu_3}\!\colon \calU(\bbZ_3) \to \!\bbQ/\bbZ$
are all non-constant.
\end{iii}\medskip

\noindent
{\bf Proof.}
i)
For~the
$\ev_{\alpha_1, \nu_p}$,
Lemma~\ref{algconst} applies, since
$d = -1$,
the linear forms
$X_0$,
$X_0+X_1$,
$X_2$,
and~$X_0+X_4$
are linearly independent modulo an arbitrary
prime~$p$,
and the cusp
$(0\!:\!0\!:\!0\!:\!1\!:\!0)$
does not satisfy the second equation
$X_0X_1 + X_2^2 = X_4X_3$
modulo any
prime~$p \neq 3$.

Similarly,~for the
$\ev_{\alpha_2, \nu_p}$,
one has
$d = -6$,
the linear forms
$2(X_0+X_2)$,
$X_0+2X_2$,
$X_1$,
and~$X_3$
are linearly independent modulo every
prime~$p$
and the cusp
$(0\!:\!0\!:\!0\!:\!0\!:\!1)$
does not satisfy the first equation
$X_0(X_0+X_1) = X_2^2 + (X_0+X_4)^2$
modulo any~prime. The assertion about
$\ev_{\alpha_1+\alpha_2, \nu_p}$~follows.\smallskip

\noindent
ii and iii)
Constancy of
$\ev_{\alpha_1, \nu_3}$
follows from Lemma~\ref{algconst}, too, since the cusp
$(0\!:\!0\!:\!0\!:\!1\!:\!0)$
is a point not lying
on~$\calU_3$.

Furthermore,~the asserted non-constancies are directly checked. In~order to do this, let us consider the points
$\smash{x = (5\!:\!\frac{12}5\!:\!-1\!:\!\frac25\!:\!1)}$
and
$\smash{x' = (-\frac{37}{13}\!:\!\frac{21}{13}\!:\!\frac4{13}\!:\!\frac5{13}\!:\!1)}$,
for which we clearly have that
$x, x' \in \calU(\bbZ_2) \cap \calU(\bbZ_3) \cap U(\bbR)$.
Explicit calculations show that
$\ev_{\alpha_1,\nu}(x) = \ev_{\alpha_2,\nu}(x) = \ev_{\alpha_1+\alpha_2,\nu}(x) = 0$
for any
place~$\nu$.
On~the other hand,
$\ev_{\alpha_1,\nu}(x') = \frac12$
for
$\nu = \nu_\infty$
and
$\nu = \nu_2$,
as well as
$\ev_{\alpha_2,\nu}(x') = \frac12$
for
$\nu = \nu_\infty$,
$\nu = \nu_2$,
and~$\nu = \nu_3$.
This~establishes ii) and~iii).\smallskip

\noindent
iv)
As~a consequence of ii) and~iii), we see that
$\ev_{\alpha_1+\alpha_2, \nu_3}$
is non-constant. Moreover, there is the point
$\smash{x'' = (-\frac{85}{91}\!:\!-\frac{15}{91}\!:\!\frac{92}{91}\!:\!\frac9{91}\!:\!1) \in \calU(\bbZ_2) \cap U(\bbR)}$,
for which direct calculations show that
$\ev_{\alpha_1+\alpha_2, \nu}(x'') = \!\frac12\,$
for
$\nu = \nu_\infty$
and
$\nu = \nu_2$.
This completes~the proof.
\eop
\end{ex}

\begin{rems}
\label{ex3_strapp}
\begin{iii}
\item
Similarly to the situation in the examples of the sections above, the Brauer classes
$\alpha_1$
and
$\alpha_2$
cause failures of strong~approximation.

There~is a violation of strong approximation off
\mbox{$S_1 := \{ p \text{~prime} \mid p \equiv 1 \pmod 4 \}$,}
due
to~$\alpha_1$.
A~\mbox{$\smash{\bbZ[\frac1{S_1}]}$-valued}
point must necessarily fulfil
$(\frac{x_0}{x_4},-1)_2 + (\frac{x_0}{x_4},-1)_\infty = 0$,
although not all adelic points outside
$S_1$
satisfy this~relation.

Similarly,~$\alpha_2$
causes a violation of strong approximation off
$$\textstyle S_2 := \{ p \text{~prime} \mid (\frac{-6}p) = 1 \} = \{ p \text{~prime} \mid p \equiv 1, 5, 7, 11 \pmodulo{24} \} \, .$$
A~\mbox{$\smash{\bbZ[\frac1{S_2}]}$-valued}
point such that
$x_0+x_2 \neq 0$
must necessarily fulfil the nontrivial relation
$\smash{(\frac{x_0+x_2}{x_4},-6)_2 + (\frac{x_0+x_2}{x_4},-6)_3 + (\frac{x_0+x_2}{x_4},-6)_\infty = 0}$.

However,~as
$U(\bbR)$
is disconnected and one of its connected components is compact, further effects are occurring. Cf.~Example~\ref{DeWe} below.

\item
There is exactly one integral point
on~$\calU$
lying on the compact component of
$U(\bbR)$,
namely~$(-1\!:\!0\!:\!1\!:\!0\!:\!1)$.

A~search for integral points 
on~$\calU$
delivered the following twelve~others.\newline
{\scriptsize
$(17\!:\!3\!:\!4\!:\!\pm13\!:\!1)$,
$(1409\!:\!147\!:\!-452\!:\!\pm383\!:\!1)$,
$(6305\!:\!12972\!:\!9043\!:\!\pm3550\!:\!1)$,
$(17741\!:\!12759\!:\!15044\!:\!\pm20351\!:\!1)$,\,%
$(-23293\!:\!-2328\!:\!-7367\!:\!\pm19622\!:\!1)$,\,%
and
$(60569\!:\!2052\!:\!11143\!:\!\pm44472\!:\!1)$.}
\end{iii}
\end{rems}

\begin{ex}[2nd continuation]
\label{DeWe}
Let\/~$U$
be as above. Then~strong approximation
off\/~$\{\infty\}$
is violated
for\/~$U$.\medskip

\noindent
{\bf Proof.}
Consider the adelic point
$\bbx$
outside
$\infty$
that is equal to
$(-1\!:\!0\!:\!1\!:\!0\!:\!1)$
at every prime
$p \neq 5, \infty$
and equal to
$(17\!:\!3\!:\!4\!:\!13\!:\!1)$
at~$p = 5$.
Working with the model
$\calU$
as above, strong approximation
off~$\{\infty\}$
would imply that there exists a sequence
$(x_n)_{n \in \bbN}$
of integral points
$x_n \in \calU(\bbZ)$
being convergent
to~$\bbx$
simultaneously with respect to the
$2$-,
$3$-,
and
\mbox{$5$-adic}
topologies. The~latter ensures that
$x_n \neq (-1\!:\!0\!:\!1\!:\!0\!:\!1)$
for~$n \gg 0$.

Moreover,
\mbox{$2$-adic}
convergence implies that
$\ev_{\alpha_1, \nu_2}(x_n) = \ev_{\alpha_1, \nu_2}(-1\!:\!0\!:\!1\!:\!0\!:\!1) = \frac12$
and
$\ev_{\alpha_2, \nu_2}(x_n) = \ev_{\alpha_2, \nu_2}(-1\!:\!0\!:\!1\!:\!0\!:\!1) = 0$
for
$n \gg 0$,
while
\mbox{$3$-adic}
convergence enforces
$\ev_{\alpha_2, \nu_3}(x_n) = \ev_{\alpha_2, \nu_3}(-1\!:\!0\!:\!1\!:\!0\!:\!1) = 0$.
The relations discussed in Remark~\ref{ex3_strapp}.i) now imply that, necessarily,
$\ev_{\alpha_1, \nu_\infty}(x_n) = \frac12$
and
$\ev_{\alpha_2, \nu_\infty}(x_n) = 0$.

In~other words,
$x_n$
must be contained in the same connected component
of~$U(\bbR)$
as the point
$(-1\!:\!0\!:\!1\!:\!0\!:\!1)$.
This~component is, however, compact and does not contain any other integral~point.
\eop
\end{ex}

\setlength\parindent{0mm}
\end{document}